\documentstyle[12pt]{article}
\catcode`\@=12
\makeatletter

\def\hexnumber@#1{\ifnum#1<10 \number#1\else
  \ifnum#1=10 A\else\ifnum#1=11 B\else\ifnum#1=12 C\else
  \ifnum#1=13 D\else\ifnum#1=14 E\else\ifnum#1=15 F\fi\fi\fi\fi\fi\fi\fi}
  \font\fivam  = msam5                
    \skewchar\fivam ='60              
  \font\fivbm  = msbm5                
    \skewchar\fivbm ='60              
  \font\sixam  = msam6                
    \skewchar\sixam ='60              
  \font\sixbm  = msbm6                
    \skewchar\sixbm ='60              
  \font\sevam  = msam7                
    \skewchar\sevam ='60              
  \font\sevbm  = msbm7                
    \skewchar\sevbm ='60              
  \font\egtam  = msam8                
    \skewchar\egtam ='60              
  \font\egtbm  = msbm8                
    \skewchar\egtbm ='60              
  \font\ninam  = msam9                
    \skewchar\ninam ='60              
  \font\ninbm  = msbm9                
    \skewchar\ninbm ='60              
  \font\tenam  = msam10               
    \skewchar\tenam ='60              
  \font\tenbm  = msbm10               
    \skewchar\tenbm ='60              
  \font\elvam  = msam10     \@halfmag 
    \skewchar\elvam ='60              
  \font\elvbm  = msbm10     \@halfmag 
    \skewchar\elvbm ='60              
  \font\twlam  = msam10   \@magscale1 
    \skewchar\twlam ='60              
  \font\twlbm  = msbm10   \@magscale1 
    \skewchar\twlbm ='60              
  \font\frtam  = msam10   \@magscale2 
    \skewchar\frtam ='60              
  \font\frtbm  = msbm10   \@magscale2 
    \skewchar\frtbm ='60              
  \font\svtam  = msam10   \@magscale3 
    \skewchar\svtam ='60              
  \font\svtbm  = msbm10   \@magscale3 
    \skewchar\svtbm ='60              
  \font\twtam  = msam10   \@magscale4 
    \skewchar\twtam ='60              
  \font\twtbm  = msbm10   \@magscale4 
    \skewchar\twtbm ='60              
  \font\twfam  = msam10   \@magscale5 
    \skewchar\twfam ='60              
  \font\twfbm  = msbm10   \@magscale5 
    \skewchar\twfbm ='60              
\newfam\msafam
\newfam\msbfam
\def\am@{\hexnumber@\msafam}
\def\bm@{\hexnumber@\msbfam}
\ifcase\@ptsize
  \textfont\msafam=\tenam  \scriptfont\msafam=\sevam
    \scriptscriptfont\msafam=\fivam
  \textfont\msbfam=\tenbm  \scriptfont\msbfam=\sevbm
    \scriptscriptfont\msbfam=\fivbm
\or
  \textfont\msafam=\elvam  \scriptfont\msafam=\egtam
    \scriptscriptfont\msafam=\sixam
  \textfont\msbfam=\elvbm  \scriptfont\msbfam=\egtbm
    \scriptscriptfont\msbfam=\sixbm
\or
  \textfont\msafam=\twlam  \scriptfont\msafam=\egtam
    \scriptscriptfont\msafam=\sixam
  \textfont\msbfam=\twlbm  \scriptfont\msbfam=\egtbm
    \scriptscriptfont\msbfam=\sixbm
\fi
\edef\yen{{\noexpand\mathhexbox\am@55}}
\edef\checkmark{{\noexpand\mathhexbox\am@58}}
\edef\circledR{{\noexpand\mathhexbox\am@72}}
\edef\maltese{{\noexpand\mathhexbox\am@7A}}
\def\ulcorner{\delimiter"4\am@70\am@70}
\def\urcorner{\delimiter"5\am@71\am@71}
\def\llcorner{\delimiter"4\am@78\am@78}
\def\lrcorner{\delimiter"5\am@79\am@79}
\mathchardef\dabar@="0\am@39
\edef\dashrightarrow{\mathrel{\dabar@\dabar@\mathchar"0\am@4B}}%

\edef\dashleftarrow{\mathrel{\mathchar"0\am@4C\dabar@\dabar@}}%
\def\Bbb{\fam\msbfam}
\mathchardef\digamma="0\bm@7A
\mathchardef\varkappa="0\bm@7B
\mathchardef\beth="0\bm@69
\mathchardef\gimel="0\bm@6A
\mathchardef\daleth="0\bm@6B
\mathchardef\hbar="0\bm@7E
\mathchardef\backprime="0\am@38
\mathchardef\hslash="0\bm@7D
\mathchardef\varnothing="0\bm@3F
\mathchardef\vartriangle="3\am@4D
\mathchardef\blacktriangle="0\am@4E
\mathchardef\triangledown="0\am@4F
\mathchardef\blacktriangledown="0\am@48
\mathchardef\square="0\am@03
\mathchardef\blacksquare="0\am@04
\mathchardef\lozenge="0\am@06
\mathchardef\blacklozenge="0\am@07
\mathchardef\circledS="0\am@73
\mathchardef\bigstar="0\am@46
\mathchardef\angle="0\am@5C
\mathchardef\sphericalangle="0\am@5E
\mathchardef\measuredangle="0\am@5D
\mathchardef\nexists="0\bm@40
\mathchardef\complement="0\am@7B
\mathchardef\mho="0\bm@66
\mathchardef\eth="0\bm@67
\mathchardef\Finv="0\bm@60
\mathchardef\diagup="3\bm@1E
\mathchardef\Game="0\bm@61
\mathchardef\diagdown="3\bm@1F
\mathchardef\Bbbk="0\bm@7C
\mathchardef\dotplus="2\am@75
\mathchardef\ltimes="2\bm@6E
\mathchardef\smallsetminus="2\bm@72
\mathchardef\rtimes="2\bm@6F
\mathchardef\Cap="2\am@65
\mathchardef\leftthreetimes="2\am@68
\mathchardef\Cup="2\am@64

\mathchardef\rightthreetimes="2\am@69

\mathchardef\barwedge="2\am@5A
\mathchardef\curlywedge="2\am@66
\mathchardef\veebar="2\am@59
\mathchardef\curlyvee="2\am@6
\mathchardef\doublebarwedge="2\am@5B
\mathchardef\boxminus="2\am@0C
\mathchardef\circleddash="2\am@7F
\mathchardef\boxtimes="2\am@02
\mathchardef\circledast="2\am@7E
\mathchardef\boxdot="2\am@00
\mathchardef\circledcirc="2\am@7D
\mathchardef\boxplus="2\am@01
\mathchardef\centerdot="2\am@05
\mathchardef\divideontimes="2\bm@3E
\mathchardef\intercal="2\am@7C
\mathchardef\leqq="3\am@35
\mathchardef\geqq="3\am@3D
\mathchardef\leqslant="3\am@36
\mathchardef\geqslant="3\am@3E
\mathchardef\eqslantless="3\am@30
\mathchardef\eqslantgtr="3\am@31
\mathchardef\lesssim="3\am@2E
\mathchardef\gtrsim="3\am@26
\mathchardef\lessapprox="3\am@2F
\mathchardef\gtrapprox="3\am@27
\mathchardef\approxeq="3\bm@75
\mathchardef\lessdot="3\bm@6C
\mathchardef\gtrdot="3\bm@6D
\mathchardef\lll="3\am@6E

\mathchardef\ggg="3\am@6F

\mathchardef\lessgtr="3\am@37
\mathchardef\gtrless="3\am@3F
\mathchardef\lesseqgtr="3\am@51
\mathchardef\gtreqless="3\am@52
\mathchardef\lesseqqgtr="3\am@53
\mathchardef\gtreqqless="3\am@54
\mathchardef\doteqdot="3\am@2B

\mathchardef\eqcirc="3\am@50
\mathchardef\risingdotseq="3\am@3A
\mathchardef\circeq="3\am@24
\mathchardef\fallingdotseq="3\am@3B
\mathchardef\triangleq="3\am@2C
\mathchardef\backsim="3\am@76
\mathchardef\thicksim="3\bm@73
\mathchardef\backsimeq="3\am@77
\mathchardef\thickapprox="3\bm@74
\mathchardef\subseteqq="3\am@6A
\mathchardef\supseteqq="3\am@6B
\mathchardef\Subset="3\am@62
\mathchardef\Supset="3\am@63
\mathchardef\sqsubset="3\am@40
\mathchardef\sqsupset="3\am@41
\mathchardef\preccurlyeq="3\am@34
\mathchardef\succcurlyeq="3\am@3C
\mathchardef\curlyeqprec="3\am@32
\mathchardef\curlyeqsucc="3\am@33
\mathchardef\precsim="3\am@2D
\mathchardef\succsim="3\am@25
\mathchardef\precapprox="3\bm@77
\mathchardef\succapprox="3\bm@76
\mathchardef\vartriangleleft="3\am@43
\mathchardef\vartriangleright="3\am@42
\mathchardef\trianglelefteq="3\am@45
\mathchardef\trianglerighteq="3\am@44
\mathchardef\vDash="3\am@0F
\mathchardef\Vdash="3\am@0D
\mathchardef\Vvdash="3\am@0E
\mathchardef\smallsmile="3\am@60
\mathchardef\shortmid="3\bm@70
\mathchardef\smallfrown="3\am@61
\mathchardef\shortparallel="3\bm@71
\mathchardef\bumpeq="3\am@6C
\mathchardef\between="3\am@47
\mathchardef\Bumpeq="3\am@6D
\mathchardef\pitchfork="3\am@74
\mathchardef\varpropto="3\am@5F
\mathchardef\backepsilon="3\bm@7F
\mathchardef\blacktriangleleft="3\am@4A
\mathchardef\blacktriangleright="3\am@49
\mathchardef\therefore="3\am@29
\mathchardef\because="3\am@2A
\mathchardef\nless="3\bm@04
\mathchardef\ngtr="3\bm@05
\mathchardef\nleq="3\bm@02
\mathchardef\ngeq="3\bm@03
\mathchardef\nleqslant="3\bm@0A
\mathchardef\ngeqslant="3\bm@0B
\mathchardef\nleqq="3\bm@14
\mathchardef\ngeqq="3\bm@15
\mathchardef\lneq="3\bm@0C
\mathchardef\gneq="3\bm@0D
\mathchardef\lneqq="3\bm@08
\mathchardef\gneqq="3\bm@09
\mathchardef\lvertneqq="3\bm@00
\mathchardef\gvertneqq="3\bm@01
\mathchardef\lnsim="3\bm@12
\mathchardef\gnsim="3\bm@13
\mathchardef\lnapprox="3\bm@1A
\mathchardef\gnapprox="3\bm@1B
\mathchardef\nprec="3\bm@06
\mathchardef\nsucc="3\bm@07
\mathchardef\npreceq="3\bm@0E
\mathchardef\nsucceq="3\bm@0F
\mathchardef\precneqq="3\bm@16
\mathchardef\succneqq="3\bm@17
\mathchardef\precnsim="3\bm@10
\mathchardef\succnsim="3\bm@11
\mathchardef\precnapprox="3\bm@18
\mathchardef\succnapprox="3\bm@19
\mathchardef\nsim="3\bm@1C
\mathchardef\ncong="3\bm@1D
\mathchardef\nshortmid="3\bm@2E
\mathchardef\nshortparallel="3\bm@2F
\mathchardef\nmid="3\bm@2D
\mathchardef\nparallel="3\bm@2C
\mathchardef\nvdash="3\bm@30
\mathchardef\nvDash="3\bm@32
\mathchardef\nVdash="3\bm@31
\mathchardef\nVDash="3\bm@33
\mathchardef\ntriangleleft="3\bm@36
\mathchardef\ntriangleright="3\bm@37
\mathchardef\ntrianglelefteq="3\bm@35
\mathchardef\ntrianglerighteq="3\bm@34
\mathchardef\nsubseteq="3\bm@2A
\mathchardef\nsupseteq="3\bm@2B
\mathchardef\nsubseteqq="3\bm@22
\mathchardef\nsupseteqq="3\bm@23
\mathchardef\subsetneq="3\bm@28
\mathchardef\supsetneq="3\bm@29
\mathchardef\varsubsetneq="3\bm@20
\mathchardef\varsupsetneq="3\bm@21
\mathchardef\subsetneqq="3\bm@24
\mathchardef\supsetneqq="3\bm@25
\mathchardef\varsubsetneqq="3\bm@26
\mathchardef\varsupsetneqq="3\bm@27
\mathchardef\leftleftarrows="3\am@12
\mathchardef\rightrightarrows="3\am@13
\mathchardef\leftrightarrows="3\am@1C
\mathchardef\rightleftarrows="3\am@1D
\mathchardef\Lleftarrow="3\am@57
\mathchardef\Rrightarrow="3\am@56
\mathchardef\twoheadleftarrow="3\am@11
\mathchardef\twoheadrightarrow="3\am@10
\mathchardef\leftarrowtail="3\am@1B
\mathchardef\rightarrowtail="3\am@1A
\mathchardef\looparrowleft="3\am@22
\mathchardef\looparrowright="3\am@23
\mathchardef\leftrightharpoons="3\am@0B
\mathchardef\leftrightharpoons="3\am@0A
\mathchardef\curvearrowleft="3\bm@78
\mathchardef\curvearrowright="3\bm@79
\mathchardef\circlearrowright="3\am@08
\mathchardef\circlearrowleft="3\am@09
\mathchardef\Lsh="3\am@1E
\mathchardef\Rsh="3\am@1F
\mathchardef\upuparrows="3\am@14
\mathchardef\downdownarrows="3\am@15
\mathchardef\upharpoonleft="3\am@18
\mathchardef\upharpoonright="3\am@16

\mathchardef\downharpoonleft="3\am@19
\mathchardef\downharpoonright="3\am@17
\mathchardef\multimap="3\am@28
\mathchardef\rightsquigarrow="3\am@20
\mathchardef\leftrightsquigarrow="3\am@21
\mathchardef\nleftarrow="3\bm@38
\mathchardef\nrightarrow="3\bm@39
\mathchardef\nLeftarrow="3\bm@3A
\mathchardef\nRightarrow="3\bm@3B
\mathchardef\nLeftrightarrow="3\bm@3C
\mathchardef\nleftrightarrow="3\bm@3D
\def\AmS{{$\cal A$\kern-.1667em\lower.5ex%
  \hbox{$\cal M$}\kern-.125em$\cal S$}}

  \font\fivfm  = eufm5                
    \skewchar\fivfm ='60              
  \font\fivfb  = eufb5                
    \skewchar\fivfb ='60              
  \font\sixfm  = eufm6                
    \skewchar\sixfm ='60              
  \font\sixfb  = eufb6                
    \skewchar\sixfb ='60              
  \font\sevfm  = eufm7                
    \skewchar\sevfm ='60              
  \font\sevfb  = eufb7                
    \skewchar\sevfb ='60              
  \font\egtfm  = eufm8                
    \skewchar\egtfm ='60              
  \font\egtfb  = eufb8                
    \skewchar\egtfb ='60              
  \font\ninfm  = eufm9                
    \skewchar\ninfm ='60              
  \font\ninfb  = eufb9                
    \skewchar\ninfb ='60              
  \font\tenfm  = eufm10               
    \skewchar\tenfm ='60              
  \font\tenfb  = eufb10               
    \skewchar\tenfb ='60              
  \font\elvfm  = eufm10     \@halfmag 
    \skewchar\elvfm ='60              
  \font\elvfb  = eufb10     \@halfmag 
    \skewchar\elvfb ='60              
  \font\twlfm  = eufm10   \@magscale1 
    \skewchar\twlfm ='60              
  \font\twlfb  = eufb10   \@magscale1 
    \skewchar\twlfb ='60              
  \font\frtfm  = eufm10   \@magscale2 
    \skewchar\frtfm ='60              
  \font\frtfb  = eufb10   \@magscale2 
    \skewchar\frtfb ='60              
  \font\svtfm  = eufm10   \@magscale3 
    \skewchar\svtfm ='60              
  \font\svtfb  = eufb10   \@magscale3 
    \skewchar\svtfb ='60              
  \font\twtfm  = eufm10   \@magscale4 
    \skewchar\twtfm ='60              
  \font\twtfb  = eufb10   \@magscale4 
    \skewchar\twtfb ='60              
  \font\twffm  = eufm10   \@magscale5 
    \skewchar\twffm ='60              
  \font\twffb  = eufb10   \@magscale5 
    \skewchar\twffb ='60              
\newfam\fmfam
\newfam\fbfam
\ifcase\@ptsize
  \textfont\fmfam=\tenfm  \scriptfont\fmfam=\sevfm
    \scriptscriptfont\fmfam=\fivfm
  \textfont\fbfam=\tenfb  \scriptfont\fbfam=\sevfb
    \scriptscriptfont\fbfam=\fivfb
\or
  \textfont\fmfam=\elvfm  \scriptfont\fmfam=\egtfm
    \scriptscriptfont\fmfam=\sixfm
  \textfont\fbfam=\elvfb  \scriptfont\fbfam=\egtfb
    \scriptscriptfont\fbfam=\sixfb
\or
  \textfont\fmfam=\twlfm  \scriptfont\fmfam=\egtfm
    \scriptscriptfont\fmfam=\sixfm
  \textfont\fbfam=\twlfb  \scriptfont\fbfam=\egtfb
    \scriptscriptfont\fbfam=\sixfb
\fi
\def\frak{\fam\fmfam}

\makeatother

\font\myaddressfont=cmti12
\font\sectionfont=cmbx12 at 12pt

 \def\lskipamount{12pt}
 \def\lskip{\vskip\lskipamount plus3pt minus2pt}
 \def\lbreak{\par \ifdim\lastskip<\lskipamount
  \removelastskip \penalty-200 \lskip \fi}

 \def\lnobreak{\par \ifdim\lastskip<\lskipamount
  \removelastskip \penalty200 \lskip \fi}

\def\section#1{\vskip 1.5truepc\centerline{\hbox {{\sectionfont #1}}}
\vskip 1truepc\noindent\stepcounter{section}}

\def\thebibliography#1{\vskip 1.5pc{\centerline {\bf References}}\vskip 4pt
\list
 {[\arabic{enumi}]}{\settowidth\labelwidth{[#1]}\leftmargin\labelwidth
 \advance\leftmargin\labelsep
 \usecounter{enumi}}
 \def\newblock{\hskip .11em plus .33em minus .07em}
 \sloppy\clubpenalty4000\widowpenalty4000
 \sfcode`\.=1000\relax}

\def\qed{{\line{\hfill $\blacksquare$}\vskip 1.5truepc}
\def\footline{\hfill}}

  \newtheorem{th}{Theorem}[subsection]
  \newtheorem{pr}[th]{Proposition}
  \newtheorem{con}[th]{Conjecture}
  \newtheorem{lem}[th]{Lemma}
  \newtheorem{df}[th]{Definition}
  \newtheorem{cor}[th]{Corollary}
  \newtheorem{dfandpr}[th]{Definition and Proposition}
 \newtheorem{rem}[th]{Remark}
 \newtheorem{ex}[th]{Example}

\def\dual{^{\ast}}
\def\N{{\Bbb N}}
\def\Zall{{{\Bbb Z}}}
\def\Z+{{{\Bbb Z}_{\geq 0}}}
\def\C{{\Bbb C}}

\def\and{{~\rm and~}}

\def\g{{\frak g}}
\def\gl{{\frak gl}}

\def\End{{{\rm End}}}
\def\dim{{{\rm dim}}}

\def\Hom{{\rm Hom}}

\def\co{^{\vee}}
\def\smallH{\bar {\cal H}}
\def\smallt{\bar {\frak t}}
\def\smallPi{\bar \Pi}

\def\smallR{\bar R}

\def\smallQ{\bar Q}

\def\smallP{\bar P}

\def\smallW{{\bar W}}
\def\smallPi{\bar \Pi}

\def\smallcalY{\bar{\cal Y}}

\def\smallS{\bar S}

\def\bigH{ {\cal H}}
\def\bigt{ {\frak t}}
\def\bigPi{ \Pi}

\def\bigR{ R}

\def\bigW{ W}
\def\bigPi{ \Pi}

\def\bigL{ {\cal L}}

\def\bigcalY{{\cal  Y}}

\def\smallS{\bar{S}}

\def\Y{Y_r}

\def\F{{\cal F}}

\def\R{R}
\def\g{{{\frak g}}}
\def\gl{{{\frak{gl}}}}
\def\C{{{\Bbb C}}}
\def\t{{{\frak t}}}
\def\e{{{\epsilon}}}
\def\N{{{\Bbb N}}}
\def\alphacheck{{\alpha}^{\vee}}
\def\echeck{\epsilon^{\vee}}
\def\T{{\cal T}}

\def\D{{\cal D}}

\def\p{\vec p}

\def\Wg{{{ W}_\g}}

\def\smalltn{{\smallt_{n}}}

\def\gln{{\frak gl}_N{(\C)}}

\def\glnhat{\widehat{\frak gl}_N(\C)}

\def\+{\oplus}
\def\*{\otimes}
\def\bigtn{{\frak t}}

\def\smalltn{\bar {\frak t}}

\def\alphacheck{\alpha^{\vee}}
\def\gen{{\hbox{\scriptsize gen}}}
\def\positive{^{\geq 0}}
\def\H{\bigH}
\def\ol{\bar}
\def\wh{\widehat}
\def\M{{\cal M}}
\def\lm{\lambda}
\def\V{{\cal V}}
\def\text{\hbox} 
\def\qed{\hfill$\square$}
\def\Z{{{\Bbb Z}}}
\begin{document}
\title{Degenerate Double Affine Hecke Algebra and Conformal Field Theory}
\author{Tomoyuki Arakawa, Takeshi Suzuki\thanks{Supported by JSPS the Research
Fellowships for Young Scientists.}\
 and
Akihiro Tsuchiya}
\date{}
\maketitle

\begin{abstract}
We introduce a class of induced representations
of the degenerate double affine Hecke algebra $\H$ of $\gl_N(\C)$
and analyze their structure mainly by means of intertwiners of $\H$ .
We also construct them from $\widehat{\frak sl}_m(\C)$-modules using 
Knizhnik-Zamolodchikov connections in the conformal field theory. 
This construction provides  natural quotients of the 
induced modules, which
correspond to the integrable 
$\widehat{\frak sl}_m(\C)$-modules.
Some conjectural formulas are presented for 
the symmetric part of them.
\end{abstract}

\begin{center}
{\bf {Introduction}}
\end{center}
In this  paper,
the representations of  the degenerate double affine Hecke
algebra ${\cal H}$
are discussed from view points of the conformal field theory
associated to the affine Lie algebra $\widehat{{\frak sl}}_m(\C)$.
The relations between the  KZ-connections of 
the conformal field theory
and the representations of  the degenerate affine Hecke algebra
was first discussed by Cherednik \cite{ch:ellip}.
And Matsuo \cite{Ma} succeeded to clarify the relations between
the differential equations satisfied by spherical functions
and KZ-connections.

At first part of this paper
we discuss the properties of
the parabolic induced modules of ${\cal H}$,
which are induced from a certain one 
dimensional representations of parabolic subalgebras of ${\cal H}$.

Secondly we will give an explicit construction of ${\cal H}$-modules
from $\widehat{\frak sl}_m(\C)$ 
modules through KZ-connections.
It will be shown that these ${\cal H}$-modules arising
from Verma modules of $\widehat{\frak sl}_m(\C)$
correspond to parabolic induced modules of ${\cal H}$.

In the final part of this paper, we will discuss the structure
of the representations of affine Hecke algebra
arising from level ${\ell}$ integrable representations of 
$\widehat{\frak sl}_m(\C)$.

We describe the contents of the paper more precisely:

In \S1 we introduce basic notions about the degenerate 
double affine Hecke algebra
$\H$.
In particular, intertwiners of weight spaces of $\H$
play essential roles in the analysis of $\H$-modules.

In \S2, we introduce parabolic induced modules  
and
investigate their structure
when the parameter is generic (Definition \ref{df;generic}), and
the results are these (Proposition \ref{pr;irr_generic},
Theorem \ref{th;generic}, Corollary \ref{cor;sympart}):

(1) The irreducibility of the standard modules are shown and 
their basis are described by intertwining operators.

(2) Decompositions of the standard modules as 
$\ol\H$-modules are obtained.

(3) The symmetric part of the standard modules are decomposed into
weight spaces with respect to the action of the center of $\ol\H$,
 and their basis are constructed again by using intertwiners.

In non-generic case, we present a sufficient condition for an induced
module
to have a unique irreducible quotient (Corollary \ref{cor;sufficient}).

\S3 is devoted to some preliminaries on affine Lie algebras and
in \S4, we realize $\H$-modules as a quotient space 
of a tensor product of 
$\g=\wh{\frak sl}_m(\C)$-modules.
More precisely,
for $\g$-modules $A,B$, we consider the space
$$\begin{array}{l}
\F(A,B)\\
=\left(A\*\*_{i=1}^N \left(\C[z_i^{\pm1}]\*\C^m\right)\* B\right)
\left/
\g'\left(A\*\*_{i=1}^N \left(\C[z_i^{\pm1}]\*\C^m\right)\* B\right) , 
\right.
\end{array}$$
where $\C[z_i^{\pm1}]\*\C^m$ is an evaluation module of $\g$ and
$\g'=[\g,\g]$
acts diagonally on the tensor product.
By combining the Knizhnik-Zamolodchikov connection with the  
Cherednik-Dunkl operator, we define an action of $\H$ on $\F(A,B)$
(Theorem \ref{th;Haction}).

In \S5, 
we construct the isomorphism between a parabolic induced module and
 $\M(\mu,\lm):=\F(M(\mu),M^*(\lm))$ for
highest and lowest Verma module $M(\mu)$ and $M^*(\lm)$,
$\lm$ and $\mu$ being weights of $\g$ (Proposition \ref{pr;standard}).

Since our construction turns out to be functorial, 
$\V(\mu,\lm):=\F(L(\mu),L^*(\lm))$ gives
a quotient module of $\M(\mu,\lm)$,
where $L(\mu)$ and $L^*(\lm)$ are the irreducible quotients 
of $M(\mu)$
and $M^*(\lm)$ respectively.
We focus on the case where $\lm$ and $\mu$ are both dominant integral
weights and
study about $\V(\mu,\lm)$.

In the last section \S6,
we focus on the symmetric part of $\V(\mu,\lm)$ for dominant integral
weights 
$\lm,\mu$,
and present a description of the basis by intertwiners 
(Conjecture \ref{con;symmetricpart})
as a consequence of the character formula $(\ref{eq;character})$,
which is still conjectural since it is proved  
under the assumption that a certain sequence of $\H$-modules 
(coming from BGG exact sequence of $\g$) is exact 
(Conjecture \ref{con;BGGexact}).   

It is interesting that character formula $(\ref{eq;character})$
also appears in the theory of the solvable lattice model \cite{ANOT}.

\section{Preliminaries }
In this section, we review the basic notions about
the degenerate double affine Hecke algebra $\bigH$.
\subsection{Affine root system }
Let $\smalltn=\bigoplus_{i=1}^N \C\echeck_i$ be 
the Cartan subalgebra of
$\gln$
with the  invariant bilinear form $(\echeck_i,\echeck_j)=\delta_{ij}$.
Define the Cartan subalgebra $\bigtn$ of the 
affine Lie algebra $\glnhat$
by
$
\bigtn=\smalltn\+\C c\+\C d
$.
Extend the non-degenerate invariant  symmetric bilinear form $(, )$ 
to $\bigtn$ by putting
$ (c,d)=1$
and $  (\e_i, c)=(\e_i,d)=(c,c)=(d,d)=0
$.
Let $\smalltn\dual=\bigoplus_{i=1}^N\C \e_i$  
be the dual space of 
$\smalltn$ 
and 
$\bigtn\dual=\smalltn \dual \oplus \C \delta \oplus \C
c^{\ast} $ 
be the dual space of $\bigtn$,
where $\e_i$, $\delta$ and  $c^{\ast} $
are   the dual vectors of $\echeck_i$,
$d$ and $c$ respectively.
We often identify 
$\bigtn\dual$ with
$\bigtn$
via  the correspondences
$ \e_i\mapsto\e_i\co$,
$\delta\mapsto c$ and $c^{\ast} \mapsto d$.
Let $\zeta\co\in\bigtn$ denote the vector corresponding to
$\zeta\in\bigtn\dual$.

Define the systems $\bigR$ of roots,  $\bigR_+$
of positive roots and  $\bigPi$ of simple roots 
of type $A_{N-1}^{(1)}$ by
$$
\begin{array}{l}
\bigR=
\left\{\alpha+k\delta\,|\,
\alpha\in\smallR,\,k\in\Zall\right\},\\
\R_{+}=\left\{\alpha+k\delta~|~
\alpha\in\smallR_{+},\,k\geq0\right\}
\sqcup \left\{-\alpha+k\delta~|~
\alpha\in\smallR_{+},\,k>0\right\},\\
\Pi=
\left\{\alpha_0:=\delta-(\e_1-\e_N)
\right\}\sqcup \smallPi,
\end{array}
$$
where $\smallR$, $\smallR_+$ and $\smallPi$ are the systems of roots,
positive roots and simple roots of type $A_{N-1}$
respectively:
$$\begin{array}{l} 
\smallR=\left\{ \alpha_{ij}=\e_i-\e_j \mid
i\ne j\right\}
,\\
\smallR_{+}=\left\{ \alpha_{ij}\mid i<j\right\},
\\ \smallPi=\left\{\alpha_1,\dots,\alpha_{N-1}\right\}
\quad (\alpha_i=\alpha_{i i+1}).
\end{array}$$

\subsection{Affine Weyl group}
Let $\bigtn'=\smalltn\+\C c\subset \bigtn$.
We  consider 
the dual space $(\bigt')\dual$ of $\bigt'$ as a subspace 
of $\bigt\dual$
via the identification
$(\bigtn')\dual=\bigtn\dual/\C \delta\cong \smalltn \dual\oplus
\C c^{\ast}$.

Let $\smallQ$ be the root lattice 
$\bigoplus_{i=1}^{N-1}\Zall \alpha_i$ and

$\smallP$ be the weight lattice $ \bigoplus_{i=1}^N \Zall\e_i$.
Let $\smallW $ be the Weyl group of $\gln$,
which is isomorphic to the symmetric group
${\frak S}_N$.
The {\em affine Weyl group} $\bigW $
is defined as a semidirect product
$$
\begin{array}{c}
\bigW =\smallW \ltimes \smallP,\\
\end{array}
$$
with the relation $w\cdot t_{ \eta}\cdot {w} ^{-1}=t_{{w}(
 \eta)}$,
where $w$ and $t_{ \eta} $
are the elements in $\bigW$ corresponding to $w\in\smallW$ and 
$\eta\in
\smallP$ respectively.
The group $\bigW$
 contains the affine Weyl group $W^a=\smallW \ltimes \smallQ$ 
of  type $A_{N-1}^{(1)}$
as its subgroup.

Let $s_{\alpha}\in\smallW$ be the reflection
corresponding to $\alpha\in\smallR $.
The action of  $\bigW $
on an element
$\xi\in\bigt$ 
is given by
the following formulas:
\begin{equation}\label{eq;ac}
\begin{array}{l}
s_{\alpha} \left(\xi\right)=
\xi-\alpha(\xi) \alpha\co \quad \left( \alpha\in\smallR  \right),\\
t_{ \eta } (\xi)=\xi+
\delta(\xi) \eta\co-
\left( \eta(\xi)+{1\over 2}( \eta, \eta)^2\delta (\xi)\right)c
\quad \left( \eta\in \smallP\right).
\end{array}
\end{equation}
The dual action on $\t^*$ is given by 
\begin{equation}\label{eq;acdual}
\begin{array}{l}
s_{\alpha}(\zeta)=
\zeta-(\alpha,\zeta) \alpha 
\quad \left( \alpha\in\smallR,\ \zeta\in\t^*  \right),\\
t_{ \eta } (\zeta)=\zeta+
(\delta,\zeta) \eta-
\left( (\eta,\zeta)+{1\over 2}( \eta, \eta)^2(\delta, \zeta)
\right)\delta
\quad \left( \eta\in \smallP,\ \zeta\in\t^*
\right).
\end{array}
\end{equation}
With respect to 
these actions, the inner products in $\t$ and $\t^*$ are $W$-invariant.
The action (\ref{eq;acdual}) preserves the set $R$ of roots.

The following action 
on $(\bigtn')^*$ is
called  the {\em affine action}:
\begin{equation}\label{eq;affineaction}
\begin{array}{l}
s_{ \alpha} \left(\zeta\right)=
\zeta-( \alpha, \zeta) \alpha \quad
\left(  \alpha\in\smallR,\zeta\in (\bigtn')\dual 
\right),\\
t_{ \eta} (\zeta)=\zeta+(\delta,\zeta)  \eta \quad \left( \eta\in 
\smallP,\zeta\in (\bigtn')\dual
\right).
\end{array}
\end{equation}


For an affine root
$\alpha=\bar\alpha+k\delta\,(\bar \alpha\in\smallR ,
k\in\Zall)$, define the corresponding affine reflection by
$s_{\alpha}=t_{-k\bar\alpha}\cdot s_{\bar\alpha}$. 
Set
$s_i=s_{\alpha_i}$
for $i=0,\dots,N-1$.
We often identify the set
$\left\{0,\dots,N-1\right\}$
with the abelian group $\Zall/N\Zall$
throughout this article.
Let $\pi=t_{\e_1}\cdot s_1\cdots s_{N-1}$. The following is 
well-known.
\begin{pr}
The group $W$ is isomorphic to the group defined by the following
generators and
relations:
$$
\begin{array}{ccl}
\hbox{generators}&:& s_i~~~(i\in\Zall/N\Zall),~\pi^{\pm 1},\\
\hbox{relations}
&:& s_i^2=1~(i\in\Zall/N\Zall),\\
&& s_i\cdot s_{i+1}\cdot s_i=
s_{i+1}\cdot s_i\cdot s_{i+1}~(i\in\Zall/N\Zall),\\
&& s_i\cdot s_{j}=s_{j}\cdot s_i~
(i-j\not\equiv\pm 1 \bmod N),\\
&& \pi\cdot s_i=s_{i+1}\cdot \pi ~(i\in\Zall/N\Zall),\\ &&
\pi\cdot\pi^{-1}=1,
\end{array}
$$
and  the subgroup $\bigW^a$ is generated by
the simple reflections $s_0,\dots,s_{N-1}$. 
In particular,
$$ 
\bigW \cong\Omega\ltimes \bigW ^a,$$
where
$\Omega=\langle \pi^{\pm 1}\rangle\cong \Zall$.
\end{pr}
For $w\in\bigW $,
let $S(w)=\bigR_{+}\cap w^{-1}(\bigR_{-})$,
where $\bigR_-$ is the set of negative roots $\bigR\backslash \bigR_+$.
The length $l(w)$ of $w\in\bigW$ is defined as the number  
$\sharp S(w)$
of the elements in $S(w)$.
For $w\in\bigW $,
an expression $w=\pi^k\cdot s_{j_1}\cdots s_{j_m}$ is called a reduced 
expression
if $m=l(w)$.
Put $\smallS(w)=S(w)\cap \smallR_{+}$.
\begin{rem}
The elements
of the
abelian subgroup $\Omega$ are characterized as elements of length $0$ 
in  $\bigW $, i.e., 
$\Omega=\left\{w\in \bigW \,|\,l(w)=0\right\}$. \end{rem}
The partial ordering $\preceq $ called the {\em Bruhat ordering} is 
defined
in the Coxeter group $\bigW^a$:
$w\preceq w'$
if  $w$ can be obtained as a subexpression
of a reduced expression  of $w'$. 
Extend this ordering $\preceq $ to
the  partial ordering  in $W$
by 
$\pi^k w\preceq \pi^{k'}w'\Leftrightarrow k= k' 
\hbox{ and } 
w\preceq w'
$
($w,w'\in W^a$).
%
\subsection{Degenerate double affine Hecke algebra}
\label{ss;DDAHA}
Let $\C[\bigW]$ denote the group algebra of $\bigW$ and
$S[\bigt']$ denote the symmetric algebra of $\bigt'$. Clearly
$\C[\bigW]=\C[\smallP]\otimes\C[\smallW]$ and 
$S[\bigt']=S[\smallt']\otimes
\C[c]$.

The degenerate double affine Hecke algebra\ was introduced by Cherednik 
(\cite{ch:ellip}). 
\begin{df}
\label{Hdef}
The degenerate double affine Hecke algebra 
(DDAHA)
$\bigH$ is the unital
associative $\C$-algebra
defined by the following properties:
\rm{(i)} As a $\C$-vector space,
$$
\bigH\cong\C[W]\otimes S[\bigt'].
$$

\noindent 
\rm{(ii)} The natural inclusions
$\C[W]\hookrightarrow \bigH$
and $S[\bigt']\hookrightarrow \bigH$
are algebra homomorphisms
(the images of $w\in\bigW$ and $\xi\in\smallt'$ will be
simply denoted by $w$ and $\xi$).

\noindent 
\rm{(iii)} The following relations hold in $\bigH${\rm :}
\begin{eqnarray}
\label{eq;rel}
& &s_i\cdot \xi-s_i
 (\xi)\cdot s_i=-\alpha_i(\xi)\quad (i=0,\dots,N-1,\,\xi\in\bigt'),
\\ 
\label{eq;rel2}
& &\pi\cdot \xi=\pi(\xi)\cdot\pi\quad (\xi\in \bigt'). 
\end{eqnarray}
\end{df}

\begin{rem}
By definition, the element $c\in\bigH$ belongs to the center ${\cal
Z}(\bigH)$.
And the original  algebra defined in 
{\rm \cite{ch:ellip}} is the
quotient algebra $\bigH(\kappa)=\bigH/\langle c-\kappa\cdot 
 \hbox{\rm{id} }\rangle$
$(\kappa \in \C\dual )$.
\end{rem}

\begin{df}
Define the {degenerate affine Hecke
algebra} {\rm(}DAHA{\rm)} $\smallH$
as the following subalgebra of $\bigH$:
$$\smallH=\langle 
w\in\smallW,\bar\xi\in\smallt \rangle \cong\C[\smallW]
\otimes S[\smallt]. $$
\end{df}

The following proposition is easy to prove (see \cite{opdam}).
\begin{pr}
\label{relation}
{\rm{(i)}}
For $w\in\bigW$ and $\xi\in\bigt'$, we have
$$\xi\cdot w=
w\left(\cdot w^{-1}(\xi)+\sum_{\alpha\in S(w)}
w(\alpha)(\xi)s_{\alpha}\right).
$$
In particular,
$$\xi\cdot w=w\cdot w^{-1}(\xi)+
\sum_{w'\prec w}c_{w'}w',$$ 
for some $c_{w'}\in\C$.

\noindent
{\rm{(ii)}} For $p\in S[\bigt']$ and $i=0,\dots,N-1$, we have
\begin{equation}
\label{eq;demazure}
p\cdot s_i-s_i\cdot s_i(p)=-\triangle_i(p),
\end{equation}
where 
$\triangle_i(p)={1\over\alpha_i\co}(p-s_i(p))\in S[\bigt']$. 
\end{pr}
($w'\in\bigW^a$).
\noindent Note that 
 $\bigH=S[\bigt']\cdot \C[\bigW]=\C [W]\cdot S[\bigt']$
from  (i) of the above proposition.
\begin{rem}
Let $S(\bigt')$ be the quotient field of $S[\bigt']$.
Then
one can naturally extend
the definition of $\bigH$ to the algebra
$\widetilde\bigH=\bigH\otimes_{S[\bigt']} S(\bigt')\cong\C[W]\otimes 
S(\bigt)$
by using {\rm Eq.\  (\ref{eq;demazure})} in place of
{\rm  Eq.\  (\ref{eq;rel})}
\end{rem}
\begin{pr}
\label{pr;center}
$\,$

\noindent
{\rm{(i)}} The center $Z(\bigH)$ of
$\bigH$
equals $\C[c]$.

\noindent
{\rm{(ii)}}
The center $Z(\smallH)$ of $\smallH$
equals 
the $\smallW$-invariant part $S[\smallt]^{\smallW}$
 of $S[\smallt]$. \end{pr}
{Outline of Proof.) \hbox{\rm{(cf. \cite{Lus:Hecke1989})}}}
(i) Let $Z_{S[\bigt']}(\bigH)$ denote the elements in $\bigH$ which
commutes
with $S[\bigt']$. By the above Proposition \ref{relation} (i),
it follows that $Z_{S[\bigt']}(\bigH)=S[\bigt']$. Then by (ii) of the 
proposition,
one can 
see that  $Z(\bigH)=S[\bigt']^{\bigW}=\C[c]$. The proof for (ii) is
 the same.\qed

\medskip
Let us identify the group algebra
$\C[\smallP]$ with the
Laurent polynomial ring
$\C[z_1^{\pm},\dots,z_N^{\pm}]$
by
putting $z_i=e^{\e_i}\in \C[\smallP]$.

The following proposition 
gives another description of the algebra $\bigH$:
\begin{pr}
\label{eq;another}
The algebra $\bigH$ is the unital associative 
$\C$-algebra such that
$$\bigH=\C[\smallP]\otimes\smallH\otimes \C[c] \quad
\text{as a }\C\text{-vector space},$$
$$
\begin{array}{l}
c\in {\cal Z}(\bigH),\\
{w}\cdot f\cdot {w}^{-1}={w}(f)~ \left(w\in \smallW,\,f\in \C
[\smallP]\right),\\ 
\left[ \bar\xi,f\right] =c\, \partial_{\bar\xi} (f) 
+\sum_{\bar\alpha\in \smallR_+}\bar\alpha(\bar\xi) 
{(1-s_{\bar\alpha})(f)
\over 1-e^{-\bar\alpha}}\cdot s_{\bar\alpha}~
\left(\bar\xi\in\smallt,\,f\in
\C[\smallP]\right), \end{array}
$$
where $\partial_{\bar\xi}\cdot e^{\bar\eta}= 
\bar\eta(\bar\xi)e^{\bar\eta}$
($\xi\in\smallt$, $\bar\eta\in\smallP$)
and
the natural inclusions
$\C[\smallP]\hookrightarrow \bigH$,
$\smallH\hookrightarrow\bigH$
and $\C[c]\hookrightarrow\bigH$
are algebra homomorphisms.
\end{pr}
{\it Proof.}  follows directly from  
Definition $\ref{Hdef}$ and Proposition \ref{relation}.
\qed

\smallskip 

Let $\bigW\positive\subset \bigW$ be the subsemigroup generated by 
$s_i (i=1,\dots,N)$
and  $\pi$. 
Define  
$$\bigH\positive=\C[\bigW\positive]\otimes S[\bigt']\cong
\C[\smallP\positive]\otimes \smallH\otimes \C[c], $$
where $\smallP\positive=\bigoplus_{i=1}^N \Z_{\geq 0} \e_i$.
Then it is easy to see that
$\bigH\positive $ is  the subalgebra of $\bigH$
by the above proposition.

\def\calC{{\cal C}}
\subsection{Intertwiners}
For an $\bigH$-module
 $V$ and $\zeta\in(\t')\dual$, define 
the weight space $V_{\zeta}$ 
and the generalized weight space $V_{\zeta}^{\gen}$
with respect to the action of $S[\bigt']$:
 $$
\begin{array}{l}
V_{\zeta}=\left\{v\in V\,|\,\xi\cdot v=\zeta(\xi)v \hbox{ for any }
 \xi\in \bigt'\right\}\\
V_{\zeta}^{\gen}=
\left\{v\in V\mid( \xi-\zeta(\xi))^k\cdot v=0\hbox{ for any }
 \xi\in \bigt',\,
\text{ for some }k \in \N\right\}.
\end{array}
$$

Let
$$
\begin{array}{l}
\varphi_i=1+s_i\alphacheck_i\in\bigH~
(i\in\left\{0,\dots,N-1\right\}\cong \Zall/N\Zall),\\ 
\varphi_{\pi}=\pi\in\bigH.
\end{array}
$$
Then
\begin{eqnarray}
\label{intint}
\varphi_i\cdot\xi=s_i(\xi)\cdot\varphi_i ,\\
\label{eq;pi}
\varphi_{\pi}\cdot \xi=\pi (\xi)\cdot \varphi_{\pi},
\end{eqnarray}
for $\xi\in\bigt'$ and $i\in\Zall/N\Zall$ (Eq. (\ref{eq;pi}) is 
nothing but the defining relation
(\ref{eq;rel2})). 
\begin{pr}$\cite{Lus:Hecke1989}$
\label{pr:int}
The above defined elements satisfy
the following relations:
$$
\begin{array}{l}
\varphi_i\cdot\varphi_{i+1}\cdot\varphi_i= 
\varphi_{i+1}\cdot\varphi_i\cdot
\varphi_{i+1}~ (i\in\Zall/N\Zall),\\
\varphi_i\cdot\varphi_{j}=\varphi_j\cdot \varphi_j~ \left(i-j\ne \pm1
\hbox{
mod }N\right),\\
\varphi_{\pi}\cdot\varphi_i=\varphi_{i+1}\cdot\varphi_{\pi}
~(i\in\Zall/N\Zall),\\
\varphi_i^2=1-{\alpha_i\co}^2~
(i\in\Zall/N\Zall).
\end{array}
$$
\end{pr}
{\it Proof.} 
We only show
$\varphi_i\cdot \varphi_{i+1}\cdot \varphi_i=\varphi_{i+1} \cdot
\varphi_i\cdot \varphi_{i+1}$.
The rest is easy.
Note that $\varphi_i$ has its inverse in $\widetilde\bigH$. 
>From Proposition \ref{relation} (i),
one can prove 
that
$\widetilde \bigH=S(\bigt')\cdot \C[\bigW]$
and
$Z_{S(\bigt')}(\widetilde\bigH):=
\left\{ X\in \widetilde \bigH\mid
\left[ X, p\right]=0 \text{ for all }p\in S(\bigt')\right\}
=S(\bigt')$. 
Since $(\varphi_i\cdot \varphi_{i+1}\cdot \varphi_i)\cdot 
(\varphi_{i+1}\cdot \varphi_{i}\cdot \varphi_{i+1})^{-1}\in
Z_{S(\bigt')}(\widetilde\bigH)$ by Eq. (\ref{intint}),
it follows that
$$\varphi_i\cdot \varphi_{i+1}\cdot \varphi_i=
p\cdot \varphi_{i+1} \cdot \varphi_i\cdot \varphi_{i+1}$$
for some $p\in S(\bigt')$.
By comparing the coefficient of $s_i \cdot s_{i+1}\cdot 
s_i=s_{i+1}\cdot
s_i\cdot s_{i+1}$
of both sides,
we get the desired equality.
\qed

Let
$$\varphi_w=\varphi_{\pi}^k\cdot
\varphi_{j_1}\cdots\varphi_{j_l}\in\bigH,$$ for 
$w=\pi^k\cdot s_{j_1}\cdots
s_{j_l}\in \bigW$ (reduced expression).
These elements are well-defined
by the Proposition
\ref{pr:int},
and
Eq.(\ref{intint})
reads as
\begin{equation}
\label{int}
\varphi_w\cdot\xi=w(\xi)\cdot\varphi_w~
\left(w\in \bigW,\,\xi\in\bigt'\right).
\end{equation}
Hence,  we have
\begin{pr}
Let $V$ be an $\bigH$-module.
Let $\zeta\in (\bigt')\dual$
and $w\in \bigW$.
Then $\varphi_w\cdot v\in V_{w(\zeta)}$
for any $v\in V_{\zeta}$.
\end{pr}
\noindent Here recall that the action of $\bigW$ 
on $(\bigt')\dual$ is the affine action (defined by Eq
(\ref{eq;affineaction})).
The element $\varphi_w$ is called the intertwiner (of weight spaces).

\begin{pr}\label{int:lead}
Let $w\in\bigW$.

\noindent {\rm (i)}
$$
\varphi_w=
w\cdot\prod_{\alpha\in S(w)}\alphacheck +
\sum_{x\prec w}x\cdot f_x  ,
$$
for some $f_x\in S[\bigt']$.

\noindent{\rm (ii)}
$$\varphi_{w^{-1}}\cdot \varphi_w=
\prod_{\alpha\in S(w)}(1-{\alphacheck}^2).$$ \end{pr}
{\it Proof.} 
(i) follows from  the well-known fact
$S(w)=\{ s_{j_l}\cdots s_{j_2}(\alpha_{j_1}),$ 
$\,s_{j_l}\cdots s_{j_3}
(\alpha_{j_2}),\dots,\alpha_{j_l}\}$ for
 $w=\pi^k s_{j_1}s_{j_{2}}\cdots
s_{j_l}$ (reduced expression).
(ii) follows from the last relation of Proposition \ref{pr:int}. \qed

\begin{lem}
\label{lem;positive}
If $w\in\bigW\positive$,
then $\varphi_w\in \bigH\positive$.
\end{lem}

 \section{Representations}
In this section we study
some important class of representations of $
\bigH$ for the next section.
\subsection{Induced representations}\label{ss;standard}
Recall that every irreducible representation of the degenerate affine
 Hecke
algebra
$\smallH$ can be obtained as the unique irreducible quotient of
a {\it standard module} (see \cite{Rogawski}),
which is  a representation induced from some
parabolic subalgebra" of $\smallH$.
Hence it is natural to start with investing
such induced modules:

Let us denote $\beta\vdash N$
if an ordered sequence $\beta=(\beta_1,\dots,\beta_m)$ 
 of positive integers 
is a (ordered) partition of N, i.e., $\sum_{i=1}^m \beta_i=N$.

For a given $\beta=(\beta_1,\dots,\beta_r)\vdash N$, 
let $I_{\beta}=( I_{\beta}^{(1)},I_{\beta}^{(2)},
\dots,I_{\beta}^{(m)})$,
where $I_{\beta}^{(k)}=\left\{
\sum_{i=1}^{k-1} \beta_{i}+1,\sum_{i=1}^{k-1} \beta_{i}+2,\dots,
\sum_{i=1}^{k} \beta_{i}\right\}$.
Define $\smallR_{\beta}=
\left\{\alpha_{i j}\in \smallR\mid
i,j\in I_{\beta_k} \text{ for some }k\right\}$,
$\smallR_{\beta.+}=\smallR_+\cap \smallR_{\beta}$
and  $\smallPi_{\beta}=\smallPi\cap \smallR_{\beta}$.
Let
$\smallW_{\beta}$ 
be the subgroup of $\smallW$ generated by
$s_{\alpha}$ ($\alpha\in\smallPi_{\beta}$).
Then clearly $\smallW_{\beta}$ is a parabolic subgroup of 
$\smallW${\rm :}
\begin{equation}
\label{parabolic}
\smallW_{\beta}\cong \smallW_{\beta_1}\times
\cdots\times\smallW_{\beta_r}
~(\smallW_{\beta_j}={\frak S}_{\beta_j}).
\end{equation}
Define  subalgebras
\begin{equation}
  \label{eq;tbeta}
\bigH_{\beta}=S[\bigt']\otimes
\C[\smallW_{\beta}]\subset \bigH,\quad
  \smallH_{\beta}=S[\smallt]\otimes\C[\smallW_{\beta}]\subset \smallH.
\end{equation}
We call $\bigH_{\beta}$ (resp.\ $\smallH_{\beta}$) 
the parabolic subalgebra of $\bigH$ (resp.\ of $\smallH$)
associated with $\beta\vdash N$.
Clearly $\bigH_{\beta}=\smallH_{\beta}\otimes \C[c]$.

Let
\begin{eqnarray}
\label{eq;t_beta}
(\bigt')\dual_{\beta}&=&
\left\{\zeta\in(\bigt')\dual\,|\,
\zeta(\alpha)=-1\hbox{ for }
\alpha\in \smallPi_{\beta}\right\}\subset (\bigt')\dual,\\
\label{eq;t_bar_beta}
\smallt\dual_{\beta}&=&
\left\{\zeta\in\smallt \dual\,|\,
\zeta(\alpha)=-1\hbox{ for }
\alpha\in \smallPi_{\beta}\right\}\subset \smallt\dual.
\end{eqnarray}
Then an element
$\zeta\in(\bigt')\dual_{\beta}$
(resp.\ $ \zeta\in\smallt\dual_{\beta}$)
defines a well-defined one-dimensional representation 
$\C { 1}_{\zeta}$ 
of 
$\bigH_{\beta}$
(resp.\ of $\smallH_{\beta}$):
$$
\begin{array}{l}
{w}\cdot {\rm 1}_{\zeta}={\rm 1}_{\zeta} 
\quad \text{for all }w\in \smallW_{\beta},\\ \xi\cdot 
{\rm 1}_{\zeta}=
\zeta(\xi){\rm 1}_{\zeta}
\quad \text{for all }\xi\in \bigt' \text{ (resp. } \xi\in \smallt \, ).
\end{array}
$$
\begin{df}
\label{df;standard}
Define the induced representation $\bigcalY_{\beta}(\zeta)$ 
of  $\bigH$ and 
 $\smallcalY_{\beta}(\zeta)$ of $\smallH$ by
$$\bigcalY_{\beta}(\zeta)=\bigH\otimes_{\bigH_{\beta}} 
\C {1}_{\zeta},
{\bar \bigcalY}_{\beta}(\zeta)=\smallH\otimes_{\smallH_{\beta}} 
\C {\bf 1}_{\zeta}.$$
The cyclic vectors $1\otimes 1_{\zeta}$ will be denoted by $\text{\bf
1}_{\zeta}$.
\end{df}
\noindent Clearly,
$$
\bigcalY_{\beta}(\zeta)\cong \C[\bigW/\smallW_{\beta}]
 \text{ as }\bigW\text{-module},\quad
\smallcalY_{\beta}(\zeta)\cong \C[\smallW/\smallW_{\beta}]
 \text{ as }\smallW\text{-module}.
$$
Let $\bar \zeta$ be the image of $\zeta\in (\bigt')_{\beta}\dual$
by the projection
 $(\bigt')_{\beta}\dual\rightarrow \smallt_{\beta}\dual$.
Then,  
$$ \bigcalY_{\beta}(\zeta)\cong
\bigH(\kappa)\otimes_{\smallH}
\smallcalY_{\beta}(\bar \zeta)\cong\C[\smallP]\otimes
\smallcalY_{\beta}(\bar \zeta),$$
where $\bigH(\kappa)=\bigH/\langle c-\kappa\, \text{id}\rangle$
and $\kappa=\zeta(c)$.
%
\subsection{Basis}
Define the following
subsets of $\bigW$ for $\beta\vdash N$:
$$
\begin{array}{l}
\bigW^{\beta}=\left\{
w\in\bigW\,|\,l(w\cdot u)\geq l(w)
\hbox{ for any }u\in \smallW_{\beta}\right\},\\
\smallW^{\beta}=\bigW^{\beta
}\cap \smallW. \end{array}$$
The following
well-known proposition will be frequently used in the rest of
this section.

\begin{pr}
\label{pr;coset}
$\,$

\noindent{\rm (i)}  
 $\bigW^{\beta}=\left\{w\in \bigW\mid S(w)\subset R_+\backslash
\smallR_{\beta,+}\right\}$ and $\smallW^{\beta}=\left\{
w\in\smallW\mid S(w)\subset \smallR_+\backslash\smallR_{\beta,+}
\right\}.$

\noindent{\rm (ii)} For any $w\in\bigW$(resp.\ $\smallW$), 
there exist a unique $w_1\in \bigW^{\beta}$(resp.\ $\smallW^{\beta}$)
and  a unique $u\in\smallW_{\beta}$, such that $w=w_1\cdot u$. 
Their length satisfy $l(w)=l(w_1)+l(u)$.
In particular, the set $\bigW^{\beta}$
(resp. $\smallW^{\beta}$)
gives a complete representatives in the coset $\bigW/\smallW_{\beta}$
(resp. $\smallW/\smallW^{\beta}$).
\end{pr}

Hence the space $\bigcalY_{\beta}(\zeta)$ has a basis $\left\{
w\cdot \text{\bf 1}_{\zeta}
\mid w\in \bigW^{\beta}\right\}$  with the partial ordering
$\preceq$
among them
induced from the Bruhat ordering.

Let us give more precise description
of this basis:
\begin{dfandpr}
\label{pr;gamma}
For $\eta\in\smallP$, define  $\gamma_{\eta}$ as the element of
$\smallW$ with shortest length possible such that 
$\gamma_{\eta}(\eta)\in\smallP_- $.
Then,
$S(\gamma_{\eta})=\left\{\alpha\in \smallR_+\,|\,
(\eta,\alpha)>0\right\}.$
\end{dfandpr}
\begin{ex}
Let $\eta=a\e_1+b(\e_2+\e_3)+c(\e_4+\e_5+\e_6)+b(\e_7+\e_8+\e_9)$
with $b<c<a$.
Then $\gamma_{\eta}=\left(
\begin{array}{ccccccccc}
1&2&3&4&5&6&7&8&9\\
9&1&2&6&7&8&3&4&5
\end{array}
\right).$
\end{ex}

\begin{lem}
\label{representative}
$\,$

\noindent{\rm (i)} For the partition $(N)$
{\rm (}the partition such that $I_{(N)}=
\left\{1,2,\dots,N\right\}${\rm )},
$$
\bigW^{(N)}=\left\{ t_\eta\cdot
\gamma_{\eta}^{-1}\mid \eta\in \smallP\right\}.
$$

\noindent{\rm (ii)}  For a general partition $\beta\vdash N$,
$$
\bigW^{\beta}=\bigW^{(N)}\cdot \smallW^{\beta}.
$$
Moreover, $l(w)=l(t_{\eta}\cdot \gamma_{\eta}^{-1})+l(x)$
for $w=t_{\eta}\cdot \gamma_{\eta}^{-1}\cdot x$
$(\eta\in \smallP$, $x\in\smallW^{\beta})$.
\end{lem}
{\it Proof.} 
follows from Proposition \ref{pr;coset} and 
Definition and Proposition \ref{pr;gamma}.
\qed
\def\al{{\alpha}}
\subsection{Representatives in  double coset}
In this subsection we want to give a description
of the double coset $\smallW\backslash \bigW/\smallW_{\beta}$,
which is related to the decomposition
of $\bigcalY_{\beta}(\zeta)$ into
a direct sum of $\smallH$-modules.
In this subsection $\beta\vdash N$ is  fixed.

Let
 $\smallP_-=\left\{
\eta\in \smallP\,|\,
(\eta,\alpha)\leq 0
\hbox{ for any }
\alpha\in\smallR_{+}
\right\}$
and 
$\smallP_{\beta,-}
=\left\{
\eta\in \smallP\,|\,
(\eta,\alpha)\leq 0
\hbox{ for any }
\alpha\in\smallR_{\beta,+}
\right\}$.

It is not difficult to see the following lemma:
\begin{lem}
\label{lem;rep_double}
The set $\smallP_{\beta,-}$ gives
 complete representatives in the double
coset $\smallW\backslash \bigW/\smallW_{\beta}$.
\end{lem}

For later purposes, let us give another description
of the representatives in the double coset: 
\begin{df}
\label{df;eta_w}
For each $w\in \smallW^{\beta}$,
define the element  $\eta_w\in \smallP_-$
by the following conditions{\rm :}
$$(\eta_w,\e_1)=0,\quad
(\eta_w,\alpha_i)=\left\{
\begin{array}{ll}
-1&\text{ if }\alpha_i\in w(\smallR_+\backslash\smallR_{\beta,+})\\
0&\text{ otherwise }
\end{array}\right.
\quad (i=1,\dots,N-1).
$$
Set $$
\begin{array}{cll}\smallP_-(w)&=&\left\{\eta+\eta_w\mid
\eta\in \smallP_-\right\}\subset \smallP_-,\\
X^{\beta}&=&
\left\{t_{\eta}\cdot w\mid w\in \smallW^\beta,\
\eta\in \smallP_-(w)\right\}.
\end{array}$$
\end{df}
Note that $\gamma_\eta=1$ for $\eta\in \smallP_-$
and thus $X^\beta\subseteq W^\beta$ 
by Lemma \ref{representative}.
\begin{lem}\label{lem;P_-(w)}
$\smallP_-(w)=
\left\{\eta\in\smallP_-\mid (\eta,\alpha)<0
\text{ for } \alpha\in w(\smallR_+\backslash \smallR_{\beta,+})
\cap \smallR_+\right\}$.
\end{lem}
{\it Proof.}
It is easy to see that the right hand side is included in
$\smallP_-(w)$.
To prove the opposite inclusion,
we suppose that $\eta\in \smallP_-(w)$ and $(\eta,\al_{ij})=0$ 
$(i<j)$
and will show that $\al_{ij}\notin  w(\smallR_+\backslash 
\smallR_{\beta,+})
\cap \smallR_+$.
Since $\eta\in\smallP_-$, 
it follows that
 $(\eta,\al_{ij})=0$ implies $(\eta,\al_k)=0$, 
and thus $\al_k\notin  w(\smallR_+\backslash \smallR_{\beta,+})$
for $i\leq k<j$.
Now $\al_{ij}\notin   w(\smallR_+\backslash \smallR_{\beta,+})\cap
\bar R_+$
follows using
$\smallR_+=w(\smallR_{\beta,+})
\sqcup \left(w(\smallR_+\backslash\smallR_{\beta,+})\cap \smallR_+
\right)
\sqcup S(w^{-1})$.
\begin{df}
\label{df;longest}
Let $\eta\in\smallP_{\beta,-}$.
We define 

$\eta_-$ $:$ the
unique element of $\smallP_-\cap \{w(\eta)\mid w\in \smallW\}$,

$\tilde\omega_\eta$ $:$ the unique longest 
element in $\smallW$ such that $\tilde\omega_\eta(\eta)=\eta_-$,

$\omega_\eta$ $:$
the unique shortest element in
$\tilde\omega_\eta \smallW_\beta$. 
\end{df}
The following lemma is easy to see.
\begin{lem}
  \label{lem;sublem}
Let $\eta\in \smallP_{\beta,-}$ and let
$w\in\smallW$ be such that $w(\eta)=\eta_-$.
Suppose that $s_i\in \smallW_\beta$ and $l(ws_i)<l(w).$
Then $s_i(\eta)=\eta$.
\end{lem}
%
\begin{lem}
\label{lem;des}
For $\eta\in \smallP_{\beta,-}$,
we have

\noindent
{\rm (i)} $\eta_-\in \smallP_-(\omega_\eta)$.

\noindent
{\rm (ii)} $\omega_\eta(\eta)=\eta_-$.

\noindent
{\rm (iii)} $\smallW t_\eta \smallW_\beta
=\smallW t_{\eta_-}\omega_\eta \smallW_\beta.$
\end{lem}
{\it Proof.} 
Put $\omega=\omega_\eta$.
To prove (i),
we shall  show
$(\eta_- -\eta_\omega,\al_i)\leq 0$ for $i=1,\dots,N-1$.
This holds obviously if $(\eta_-,\al_i)<0$
or $(\eta_\omega,\al_i)=0$.
We suppose that 
$$(\eta_-,\al_i)=0\ \text{and }
\al_i\in \omega(\smallR_+),$$
and we shall show 
$\al_i\in \omega(\smallR_{\beta,+}).$
Note that this implies $
(\eta_-,\al_i)=0 \Rightarrow (\eta_\omega,\al_i)=0$
and thus completes the proof of $\eta_-\in \smallP_-(\omega)$.
Write $\omega=\tilde\omega z$, where 
$\tilde\omega=\tilde\omega_\eta$ and
$z\in \smallW_\beta$.
Since $(\eta_-,\al_i)=0$ we have 
$s_i\tilde\omega(\eta)=s_i(\eta_-)=\eta_-$.
Hence $l(s_i\tilde\omega)<l(\tilde\omega)$, and
this is equivalent to
$\tilde\omega^{-1}(\al_i)\in \smallR_-.$
Combining with $\al_i\in\omega(\smallR_+)$,
we have
 $\omega^{-1}(\al_i)\in \smallR_+\cap z^{-1}(\smallR_-)$,
which is a subset of
$\smallR_{\beta,+}$ as $z\in\smallW_\beta$.

Let us prove (ii).
Since 
$\omega(\eta)=\tilde\omega z(\eta)$,
it is enough to prove $z(\eta)=\eta$.
If $z=s_{i_1}s_{i_2}\cdots s_{i_p}$ 
is a reduced expression of $z$,
then we have $s_{i_1},s_{i_2},
\dots,s_{i_p}\in \smallW_\beta$
and $l(\tilde\omega)>l(\tilde\omega s_{i_1})
> l(\tilde\omega s_{i_1}s_{i_2})>\cdots 
> l(\tilde\omega s_{i_1}\cdots s_{i_p})=l(\tilde\omega z)$. 
Now the statement follows from Lemma \ref{lem;sublem}.

The statement (iii) is a direct consequence of (ii).
\qed
\begin{pr}
\label{pr;representative2}
The set $X^{\beta}\subseteq \bigW$ 
gives complete representatives in the double coset
$\smallW\backslash \bigW/\smallW_{\beta}$.
\end{pr}
{\it Proof.} 
By Lemma \ref{lem;des},
we can define
the map $f:\smallP_{\beta,-}\to X^\beta$ by
$f(\eta)=t_{\eta_-}\omega_\eta$ 
($\eta\in \smallP_{\beta,-}$).
We define the map $g:X^\beta\to \smallP$
by 
$g(t_\eta w)= w^{-1}(\eta)$ 
($w\in \smallW^\beta$,
$\eta\in \smallP_{-}(w)$).
Let $\al\in\smallR_{\beta,+}$. 
Then $w(\al)\in \smallR_+$ by Proposition
\ref{pr;coset}-(i).
Therefore
$(w^{-1}(\eta),\al)=(\eta,w(\al))\leq 0$ and thus
$g(t_\eta w)\in \smallP_{\beta,-}$.
Namely the image of $g$ is included 
in $\smallP_{\beta,-}$.
By Lemma \ref{lem;des}-(ii), 
we have $g\circ f={\rm{id}}_{\smallP_{\beta,-}}$.
Let us check $f\circ g={\rm{id}}_{X^\beta}$.
For $t_{\eta'} w\in X^\beta$,
put $\eta=w^{-1}(\eta')\in\smallP_{\beta,-}$.
Then $f\circ g(t_{\eta'} w)=t_{\eta_-}\omega_{\eta}$.
It is easy to see $\eta'=\eta_-$.
We will show $w=\omega_{\eta}.$
Since $\eta'\in\smallP_-(w)$,
it can be shown using Lemma \ref{lem;P_-(w)}
that
$$S(w)=\{ 
\al\in\smallR_+\setminus\smallR_{\beta,+}\mid
(\eta,\al)\geq 0\}.
$$
Similarly, by $\eta_-=\eta'\in \smallP_-(\omega_\eta)$,
we have 
\begin{equation}\label{eq;formula_longest}
S(\omega_{\eta})
=\{\al\in\smallR_+\setminus\smallR_{\beta,+}\mid
(\eta,\al)\geq 0\}=S(w).
\end{equation}
Thus we have $w=\omega_\eta$ as required.
Now the statement follows from Lemma \ref{lem;rep_double}
and Lemma \ref{lem;des}-(iii).
\qed
\begin{rem}\label{rem;another}
Let $\eta\in \smallP_{\beta,-}$. The formula
$(\ref{eq;formula_longest})$ implies that
$\omega_\eta$ is also characterized as
the longest element such that
\begin{equation}
\omega_\eta\in\smallW^\beta\text{ and }
  \omega_\eta(\eta)=\eta_-.
\end{equation}
\end{rem}
\def\one{{{\text{\bf 1}}}}
\begin{df}
For $w\in W$,
define $\beta(w)$ as an ordered partition
of $N$ such that
$$\smallPi_{\beta(w)}=\smallPi\cap w(\smallR_{\beta,+}).$$
Set
$$\begin{array}{cll}
 \bar W_{\beta(w)}&=&\langle s_\al\mid \al\in \bar\Pi_{\beta(w)}\rangle,\\
\smallW_{[w]}&=&\{u\in\smallW\mid uw\smallW_\beta=w\smallW_\beta\}.
\end{array}$$ 
\end{df}
The proof of the following lemma is elementary.
\begin{lem}\label{lem;Pibetax}
{\rm (i)}
For $w\in W^\beta$, the followings are equivalent$:$
$$(a)\ s_iw\notin W^\beta.
\quad
(b)\ s_i\in \bar W_{\beta(w)}.
\quad
(c)\ s_i\in \bar W_{[w]}.$$

\noindent
{\rm (ii)} Let $x=t_\eta w\in X^\beta$
$(\eta\in \bar P_-(w)$, $w\in\bar W^\beta)$. 
Then
$$\bar W_{[x]}=\bar W_{\eta}\cap \bar W_{[w]},$$
where $\bar W_{\eta}=\{u\in \bar W\mid w(\eta)=\eta\}.$
\end{lem}
\begin{pr}
  \label{pr;a}
For $x\in X^\beta$, we have
$$\smallW_{[x]}=\smallW_{\beta(x)}.$$
\end{pr}
{\it Proof.}
Let $x=t_\eta w\in X^\beta$ $(\eta\in \bar P_-(w),\ w\in\bar W^\beta)$.
It is obvious that
$\smallW_{\beta(x)}\subseteq\smallW_{[x]}$.
We will prove $u\in\smallW_{[x]}
\Rightarrow u\in\smallW_{\beta(x)}$
 by the induction on $l(u)$.
The case $l(u)=1$ has been proved in Lemma
\ref{lem;Pibetax}-(i).
Let $l(u)>1$.
It is enough to show that there exists $s_i\in \smallW_{[x]}$ 
such that  $s_iu\in  \smallW_{[x]}$ and
$l(u)=l(s_iu)+1$.
Because $\bar W_\eta$ is 
generated by the simple refrections
$s_i$ such that $(\eta,\al_i)=0$,
we can find $s_i\in \bar W_{\eta}$
such that $l(u)=l(s_i u)+1$.
We will show 
$s_i\in \bar W_{[x]}$, or equivalently,
$v:=s_iu\in \bar W_{[x]}$.

If $s_i x\notin W^\beta$, then $s_i\in \bar W_{[x]}$
by Lemma \ref{lem;Pibetax}-(i).
Suppose  $s_i x\in W^\beta$.  
Then, noting
\begin{equation}
\label{eq;R_+}
\bar R_+=w(\bar R_{\beta,+})
\sqcup (w(\bar R_+\setminus \bar R_{\beta,+})\cap \bar R_+)
\sqcup S(w^{-1}),
\end{equation}
we have $s_i\in S(w^{-1})$.
It follows from $u=s_iv\in\bar W_{[x]}
\subseteq \bar W_{[w]}$ that $S(w^{-1})\subseteq
S(w^{-1}v^{-1}s_i)$ and in particular,
$w^{-1}v^{-1}s_i(\al_i)=-w^{-1}v^{-1}(\al_i)\in \bar R_-$.
Hence 
$v^{-1}(\al_i)\notin S(w^{-1}).$
Since $(\eta,v^{-1}(\al_i))=(v(\eta),\al_i)=0$,
we have 
 $v^{-1}(\al_i)\in w(\bar R_{\beta,+})$
(by using (\ref{eq;R_+}) again).
Hence $w\bar W_\beta=uw\bar W_\beta=s_ivw\bar W_\beta
=vs_{v^{-1}(\al_i)}w\bar W_\beta=vw\bar W_\beta$.
Therefore $v\in \bar W_{[x]}$ by Lemma \ref{lem;Pibetax}-(ii).
\qed
\begin{cor}\label{cor;rep2}
Let $x\in X^\beta$.
We have $\smallW^{\beta(x)} x\subseteq W^\beta$.
Moreover,
the map 
$$\bigsqcup_{x\in X^\beta} \bar{W}^{\beta(x)}
\to W^\beta$$ 
given by
$w\mapsto wx$ $(w\in \bar{W}^{\beta(x)})$ is bijective.
\end{cor}
\subsection{Structure in generic}
In this section $\beta\vdash N$ is 
still fixed.
\begin{df}\label{df;generic}
A weight $\zeta\in (\bigt')\dual_{\beta}$ {\rm (}resp.\
$\smallt_{\beta}\dual${\rm )} is said to be $\beta$-generic
if
$(\zeta,\alpha)\not\in \left\{1, 0,-1\right\}$
for any $\alpha\in \cup_{w\in \bigW^{\beta}} S(w)$
{\rm (}resp.\ $\alpha\in \cup_{w\in \smallW^{\beta}} \bar S(w)${\rm )}.
\end{df}
\begin{rem}
A weight
$\zeta\in(\bigt')_{\beta}\dual$ 
$($resp.\ $ \smallt_{\beta}\dual)$
 is $\beta$-generic if and only if
$
(\zeta, \alpha)\not\in
\left\{1,0,-1\right\}
$ for any $ 
\alpha\in R_+\backslash\smallR_{\beta,+}
$ 
$($resp.\ $ \smallR_+\backslash\smallR_{\beta,+})$.
\end{rem}
\begin{pr}
\label{pr;irr_generic}
If $\zeta\in (\bigt')\dual_{\beta}$
{\rm (}resp.\
$\smallt_{\beta}\dual${\rm )}
is $\beta$-generic,
then $\bigcalY_{\beta}(\zeta)$ 
{\rm (} resp.\
$\smallcalY_{\beta}(\zeta)${\rm )}
is an irreducible $\bigH$-module
{\rm (} resp.\
$\smallH$-module{\rm )}
with a basis 
$\left\{ 
\varphi_w\cdot  {\bf 1}_{\zeta}\,
|\, w\in\bigW^{\beta} \right\}$
{\rm (resp.\ }
$\left\{ 
\varphi_w\cdot {\bf 1}_{\zeta}\,
|\, w\in\smallW^{\beta} \right\}${\rm )}.
\end{pr}
{\it Proof.} 
Let $\zeta\in (\bigt')\dual_{\beta}$.
By  Proposition \ref{int:lead} (i) and Proposition \ref{pr;coset} (ii),
each vector $\varphi_w\cdot \text{\bf 1}_{\zeta}$ has the  form of 
$\varphi_w \text{\bf 1}_{\zeta}=
c\cdot w\text{\bf 1}_{\zeta}+\sum_{v\prec w}
c_{v}v \text{\bf 1}_{\zeta}$ for some $c\ne 0$ and $c_v\in \C$.
Hence 
the elements  
$\left\{\varphi_w\cdot \text{\bf 1}_{\zeta}\mid w\in
\bigW^{\beta}\right\}$
form a basis of $\bigcalY_{\beta}(\zeta)$.
Moreover, by Proposition \ref{int:lead} (ii),
each $\varphi_w$ ($w\in\bigW^{\beta}$) is invertible on
the vector $\text{\bf 1}_{\zeta}$.
The proof of the statement about
 $\smallcalY_{\beta}(\zeta)$ is the same.
\qed

\medskip
Any weight $\zeta\in (\t')^*$ defines a
character
$[\zeta]:S[\t']^{\bar W}=S[\bar\t]^{\bar W}\otimes\C [c]\to\C$
via the evaluation at $\zeta$.
It is regarded as an element of $(\t')^*/ \bar W$.
(Recall that $S[\bar\t]^{\bar W}$ is the center of $\bar \H$.)
\begin{pr}
  Suppose that $\zeta\in (\t')_\beta^*$ $($resp. $\bar\t_\beta^*)$
is $\beta$-generic
and $(\zeta,\delta)\ne 0$.
Then

\noindent
{\rm{(i)}} 
The map
$W\to (\t')^*$ $($resp 
$\bar W\to  \bar\t^*)$
given by $w\mapsto w(\zeta)$ is injective.

\noindent{\rm{(ii)}}
The map $X^\beta\to (\t')^*/\bar W$
given by $x\mapsto [{x(\zeta)}]$ is injective.
\label{pr;weight_inj}
\end{pr}
{\it Proof.}
The $\beta$-genericity of $\zeta\in
(\t')_\beta^*$ implies that
$(\zeta,\al)\neq 0$ for any $\al\in R$.
>From this, the statements follow easily.
\qed
%
\begin{lem}\label{lem;Wbetax}
Let $\zeta\in(\t')^*_\beta$. 
Let $w\in W^\beta$ and
$\al_i\in \smallPi_{\beta(w)}$.
Then
$s_i\cdot\varphi_w\cdot\one_\zeta=\varphi_w\cdot\one_\zeta$
for the cyclic vector
 $\one_\zeta
\in{\cal Y}_\beta(\zeta)_\zeta$.
\end{lem}
{\it Proof.}
Since
 $s_i w=w s_{w^{-1}(\al_i)}$ and 
 ${w^{-1}(\al_i)}\in\smallR_{\beta,+}$,
we have $l(w s_{w^{-1}(\al_i)})=l(s_i w)=l(w)+1$.
By Proposition
\ref{pr;coset}-(ii), we have $l(w s_{w^{-1}(\al_i)})
=l(w)+l(s_{w^{-1}(\al_i)}).$
Hence
$\varphi_i\cdot \varphi_w \cdot\text{\bf 1}_{\zeta}
=\varphi_{s_i w} \cdot\text{\bf 1}_{\zeta}
=\varphi_{w s_{w^{-1}(\al_i)}}\cdot\text{\bf 1}_{\zeta}
=\varphi_w\cdot \varphi_{s_{w^{-1}(\alpha_i)}}
 \cdot\text{\bf 1}_{\zeta}.$
Noting $({w(\zeta)},\al_i)=-1$ and $s_{w^{-1}(\al_i)}
 \text{\bf 1}_{\zeta}= \text{\bf 1}_{\zeta}$,
we have
$$(1-s_i)\varphi_w \text{\bf 1}_{\zeta}
=\varphi_i\cdot \varphi_w \text{\bf 1}_{\zeta}
= \varphi_w \cdot \varphi_{s_{w^{-1}(\alpha_i)}} 
\text{\bf 1}_{\zeta}=
\varphi_w\cdot (1-s_{w^{-1}(\alpha_i)})
\text{\bf 1}_{\zeta}=0$$ as required.
\qed
\begin{th}
\label{th;generic}
Suppose that $\zeta
\in(\bigt_{\beta}')\dual$
is $\beta$-generic and $(\zeta,\delta)\neq 0$.
Then
\begin{equation}\label{eq;decomposition_generic}
\bigcalY_{\beta}(\zeta)=
\bigoplus_{x\in X^{\beta}}\bar \H \cdot\varphi_x\one_\zeta.
\end{equation}
Moreover, each 
$\bar \H \cdot\varphi_x\one_\zeta$
is isomorphic to
$\smallcalY_{\beta(x)}(\bar{x(\zeta)})$
and irreducible as an $\smallH$-module.
\end{th}
{\it Proof.} 
First we will prove the latter part. 
It is easy to check that
$\bar{x(\zeta)}\in\bar\t_{\beta(x)}^*$ 
and it 
is $\beta(x)$-generic.
Hence $\smallcalY_{\beta(x)}(\bar{x(\zeta)})$
is irreducible. 
By Lemma \ref{lem;Wbetax},
there exists a surjective $\bar{\H}$-homomorphism
$$\bar{\cal Y}_{\beta(x)}(\bar{x(\zeta)})
\to \bar{\H}\varphi_x\one_\zeta$$
such that $\one_{\bar{x(\zeta)}}\mapsto 
\varphi_x\one_\zeta$.
This is bijective since 
$\bar{\cal Y}_{\beta(x)}(\bar{x(\zeta)})$ is irreducible.

It follows from 
Proposition \ref{pr;irr_generic} and Proposition \ref{pr;weight_inj}-(i)
that $\bar{\H}\varphi_x\one_\zeta=\oplus_{w\in \bar{W}^{\beta(x)}}
\C\varphi_{wx}\one_{\zeta}$.
Now the statement follows from Corollary \ref{cor;rep2}
and Proposition \ref{pr;irr_generic}.
\qed

\medskip
For an $\H$-module $V$, define its $\bar W$-invariant part
by $V^{\bar W}=\{v\in V\mid wv=v\hbox{ for all }w\in \bar W\}$.
Then 
$S[\t']^{\bar W}$ acts on $V^{\bar W}.$
For a character $\chi:S[\t']^{\bar W}\to\C$, we set
$$V^{\bar W}_{\chi}=\{
v\in V^{\bar W}\mid p\cdot v=\chi(p) v\hbox{
 for all }p\in S[\t']\}.$$
By Proposition \ref{pr;weight_inj}-(ii) and Theorem 
\ref{th;generic}, we have
\begin{cor}
\label{cor;sympart}
Suppose that
 $\zeta\in (\t')_\beta^*$ 
is $\beta$-generic and $(\zeta,\delta)\ne 0$.
Then
$$
\begin{array}{rl}
\bigcalY_{\beta}(\zeta)^{\smallW}
=&\bigoplus_{x\in X^\beta}
\bigcalY_{\beta}(\zeta)_{ 
[x(\zeta)]}^{\smallW}
\end{array}
$$
Moreover, for each $x\in X^\beta$, we have
$\bigcalY_{\beta}(\zeta)_{ 
[x(\zeta)]}^{\smallW}
=\C Q\varphi_x\one_{\zeta}\ne 0$,
where  $Q={1\over N!}\sum_{w\in\smallW}w\in\C[\bar W]$.
\end{cor}
\subsection{Unique irreducible quotients}
\begin{lem}
\label{lem;decom_gen}
Let $\beta\vdash  N$ and $\zeta\in(\bigt')\dual_{\beta}$.
Then,
$$\bigcalY_{\beta}(\zeta)
=\bigoplus_{\zeta'\in\left\{w(\zeta)\mid w\in\bigW^{\beta}\right\}}
\bigcalY_{\beta}(\zeta)_{\zeta'}^{\hbox{\rm \gen}}$$
and {\rm dim} $\bigcalY_{\beta}(\zeta)_{\zeta'}^{\hbox{\rm \gen}}=
\sharp\left\{w\in\bigW^{\beta}\mid w(\zeta)=\zeta'\right\}$
if it is finite.
\end{lem}
{\it Proof.} 
Use the fact
 from
 Proposition \ref{relation} (i),
 that for any $\xi\in \bigt'$ and $w\in\bigW^{\beta}$,
$\xi\cdot w\text{\bf 1}_{\zeta}=
w(\zeta)(\xi)w{\bf 1}_{\zeta}+
\sum_{v\prec w}
d_{v}v$ for some $d_v\in \C$.
\qed

\begin{pr}
\label{pr;unique_quotient}
Let $\beta\vdash  N$ and $\zeta\in(\bigt')\dual_{\beta}$.
Suppose that {\rm dim} $\bigcalY_{\beta}(\zeta)_{\zeta}=1$.
Then the $\bigH$-module $\bigcalY_{\beta}(\zeta)$ has an unique 
irreducible quotient $\bigL_{\beta}(\zeta)$.
\end{pr}
{\it Proof.} 
Let ${ M}$ be a proper submodule of 
$\bigcalY_{\beta}(\zeta)$.
Then $M$ admits the (generalized) weight decomposition,
thus
${ M}\subset \bigoplus\limits_{
\zeta'\ne \zeta}
\bigcalY_{\beta}(\zeta)_{\zeta'}^{\hbox{\rm \gen}}$
since
the vector $\text{\bf 1}_{\zeta}$ is a cyclic vector
and if $M_{\zeta}^{\hbox{\rm \gen}}\ne \{0\}$,
then $M_{\zeta}\ne\{0\}$.
Hence the sum of the all proper submodules 
 is the maximal proper submodule of $\bigcalY_{\beta}(\zeta)$.
\qed
\def\Z{{{\Bbb Z}}}

\begin{pr}
\label{pr;w(zeta)=zeta}
Let $\beta=(\beta_1,\dots,\beta_m)\vdash  N$ 
and $\zeta\in(\bigt')\dual_{\beta}$.
Set $k_a=\sum_{i=1}^{a-1}\beta_i+1$
$(a=1,\dots,m)$. 
Suppose that the following  conditions hold{\rm :}
\begin{equation}
\label{cond;dominant}
\begin{array}{l}
(\zeta,\delta)\notin {\Bbb Q}_{\leq 0},\\
(\zeta,r\delta+\alpha_{k_a\,k_{b}})\notin\Z_{\leq 0}
\hbox{ {\rm  for any } } 1\leq a<b\leq m,\
 r\in\Z_{\geq 0},\\ 
(\zeta,r\delta-\alpha_{k_a\,k_b})\notin \Z_{\leq 0}
\hbox{ {\rm  for any } } 1\leq a<b\leq m,\
 r\in\Z_{>0}.
\end{array}
\end{equation}
Then 
$
w(\zeta)\ne \zeta 
$ for any $w\in\bigW^{\beta}\setminus\{1\}$,
where $1$ denotes the unit element of $W$.
In particular ${\cal Y}_\beta(\zeta)_\zeta^{\hbox{{\rm \gen}}}
=\C\one_\zeta$.
\end{pr}
{\it Proof.} 
(Step 1)
First we prove that 
$w(\zeta)\ne \zeta $ for any 
$w\in\smallW^{\beta}\backslash \{ 1\}.$
Suppose that 
$w(\zeta)=\zeta$ for $w\in \smallW^\beta\setminus \{1\}$.
Let $p$ be the smallest number such that $w(p)\ne p$,
and let $k_a\leq p < k_{a+1}$.
Then
 $w^{-1}(p)=k_b>p(\geq k_a)$ for some $b$ since $w\in\smallW^{\beta}$.
But then $0=(\zeta-w(\zeta),\e_p)=(\zeta,\alpha_{p\, w^{-1}(p)})
=(\zeta,\al_{p\, k_a})+(\zeta,\al_{k_a,k_b})$.
This implies $(\zeta,\al_{k_a,k_b})=
(\zeta,\al_{k_a,p})\in \Z_{\leq 0}$.
This contradicts the condition (\ref{cond;dominant}).

(Step 2)
Let $x\in W^\beta$.
By Lemma \ref{representative}, we can write
$x=t_{\eta}\cdot\gamma_{\eta}^{-1}
\cdot w= \gamma_{\eta}^{-1} \cdot t_{\gamma_{\eta}(\eta)} \cdot w
$
$(\eta\in\smallP,\,w\in\smallW^{\beta})$.
Suppose that $x(\zeta)=\zeta$.
By (Step 1), it is sufficent to prove 
$\eta = 0 $. 
Suppose $\eta \neq 0 $.
Putting $\e=\sum_{i=1}^N\e_i$, 
we have
$
(\zeta,\e)=
(x(\zeta),\e)=(t_{\gamma_{\eta}(\eta)}w(\zeta),\e)
=(\zeta,\e)+(\zeta,\delta)
({\gamma_{\eta}(\eta)},\e),
$
and thus $({\gamma_{\eta}(\eta)},\e)=0$.
This implies  $r:=-(\gamma_{\eta}(\eta),\e_1)\in \Z_{> 0}$ as
$\gamma_{\eta}(\eta)\in\smallP_-$.
Since 
$t_{\gamma_{\eta}(\eta)} w(\zeta)=\gamma_{\eta}(\zeta)$,
we have
\begin{equation}\begin{array}{rl}
0&=(t_{\gamma_{\eta}(\eta)} w(\zeta)- \gamma_{\eta}(\zeta),\e_1)
=(\gamma_{\eta}(\eta),\e_1)(\zeta,\delta)+
(\zeta,\alpha_{w^{-1}(1),\gamma_{\eta}^{-1}(1)}).
\end{array}\label{eq;contradiction}
\end{equation}
Note that 
$w^{-1}(1)=k_a$ for some 
$a$.
Write $\gamma_\eta^{-1}(1)=p$ and
let $k_b\leq p< k_{b+1}$.
Then  (\ref{eq;contradiction})
leads $0=-(\zeta, r\delta-\al_{k_a,p})$
and thus
\begin{equation}
\label{eq;neg}
(\zeta,r\delta-\al_{k_a,k_b})=(\zeta,\al_{k_b,p})\in\Z_{\leq 0}.
\end{equation}
It is easy to see that (\ref{eq;neg}) never occurs
under
the condition (\ref{cond;dominant}).
\qed

\medskip
As a consequence of the results above, we have
\begin{cor}
\label{cor;sufficient}
Let $\beta\vdash N$. Suppose that $\zeta\in(\bigt')_{\beta}\dual$
satisfies the condition $(\ref{cond;dominant})$.
Then the $\bigH$-module $\bigcalY_{\beta}(\zeta)$ has a unique
irreducible quotient $\bigL_{\beta}(\zeta)$.
\end{cor}
\def\b{{\frak b}}
\def\e{{\epsilon}}
\def\g{{\frak g}}
\def\h{{\frak h}}
\def\j{{\jmath}}
\def\k{{\kappa}}
\def\l{{\ell}}
\def\m{{m}}
\def\n{{\frak n}}
\def\p{{\vec{p}}}
\def\r{{\frak r}}
\def\s{{\sigma}}
\def\t{{\frak t}}
\def\v{{v}}
\def\y{{p}}
\def\A{{a}}
\def\C{{\Bbb C}}
\def\D{{\cal D}}
\def\F{{\cal F}}
\def\Sym{{\frak S}}
\def\H{{\cal H}}
\def\L{{\cal L}}
\def\M{{\cal M}}
\def\O{{\cal O}}
\def\P{{\cal P}}
\def\Q{{\Bbb Q}}
\def\T{{\cal T}}
\def\V{{\cal V}}
\def\W{{\frak W}}
\def\Y{{\cal Y}}
\def\Z{{\Bbb Z}}
\def\+{\mathop{\oplus}}
\def\*{\mathop{\otimes}} 
\def\SumN{\mathop{\sum}_{i=1}^N}
\def\Sum+{\mathop{\sum}_{\alpha\in R_+}}
\def\dsumN{\mathop{\oplus}_{i=1}^N}
\def\ol{\bar}
\def\ij{{\act_{ij}}}
\def\ji{{\act_{ji}}}
\def\alg{{\gamma}}
\def\pair{{\lm,\mu}}
\def\clpair{{\bar\lm,\bar\mu}}
\def\one{{\text{{\rm \bf 1}}}}
\def\act{{\theta}}
\def\trans{{T}}
\def\smallm{{{\scriptscriptstyle m}}}
\def\proj{{\varpi}}
\def\car{{\frak k}}
\def\rhog{{{ {\rho}\ch{} }_{\bar\g}}}
\def\ba{{{b}}}
\def\OX{{\cal R}}
\def\Og{{\cal O}}
\def\Rg{{R_{\ol\g}}}
\def\Pg{{P_\g}}
\def\CP{\C[P]}
\def\Hh{\wh{\H}}
\def\gl{{{{\frak {gl}}}}}
\def\SL{{{\frak sl}}}
\def\Vbox{V_{\Box}}
\def\VN{{\Vbox^{\* N}}}
\def\KZ{\nabla}
\def\Wg{{W_{\frak g}}}
\def\Lm{\Lambda}
\def\lm{\lambda}

\def\al{\alpha}
\def\ch{\,\check{}}
\def\wh{\widehat}
\def\bra{{\langle}}
\def\ket{{\rangle}}
\def\ge0{{\geq 0}}
\def\End{{\text{{\rm End}}}\,}
\def\Hom{{\text{{\rm Hom}}}\,}
\def\dim{{\text{{\rm dim}}}\,}
\def\char{{\text{{\rm ch}}}\,}
\def\wt{{\text{{\rm wt}}}\,}
\section{Affine Lie algebras and classical $r$-matrices}
We introduce some notations and basic facts on affine Lie algebras,
which will be used in the next section. 
\subsection{Affine Lie algebra of type $A_{\m-1}^{(1)}$}
%
%
Let $\ol\g$ be the Lie algebra $\SL_\m(\C)$ 
and $\g=\ol\g\*\C[t,t^{-1}]\+\C c_\g\+\C d_\g$ be
the affine Lie algebra 
$\wh\SL_\m(\C)$ associated with $\ol\g$
with the commutation relations
$$
\begin{array}{l}
  [X\* f,Y\* g]=[X,Y]\* fg+(X\* f,Y\* g)_\g c_\g,  \\
 \left[d_\g,X\* f\right]=X\*t{df\over dt},\quad c_\g\in Z(\g)
\end{array}
$$
for $X,Y\in\ol\g$ and $f,g\in\C[t,t^{-1}]$, where 
the invariant bilinear form $(\ ,\ )_\g$ is defined by
\begin{equation}
\begin{array}{c}
   (X\*f,Y\*g)_\g=\hbox{tr}(XY)
   \mathop{\hbox{Res}}_{t=0} f\frac{dg}{dt}dt,\\
   (X\*f,c_\g)_\g=(X\*f,d_\g)_\g=0,\\
   (c_\g,d_\g)_\g=1,\ (c_\g,c_\g)_\g=(d_\g,d_\g)_\g=0.
  \end{array}
\end{equation}
A Cartan subalgebra $\h$ of $\g$
 is given by
 $ \h=\ol\h \+\C c_\g\+\C d_\g$,
where $ \ol\h$ is a Carten subalgebra of $\ol\g$.
Its dual space is denoted by $\h^*=\bar\h^*\+\C c^*_\g\+\C\delta_\g$,
where $c^*_\g$ and $\delta_g$ denote the dual of
$d_\g$ and $c_\g$ respectively:
 $c^*_\g(x)=(d_\g,x)_\g,\ \delta_g(x)=(c_\g,x)_\g
$ for any $x\in\h$.
We write the root system and the set of positive roots
of $\ol\g$ by $\Rg$ and $\Rg_+$ respectively.
Then we have the root space decomposition  
$  \ol\g=\ol\h\+(\+_{\alg\in\Rg}\ol\g_\alg)$.
We choose the set 
$\{E_\alg\in\ol\g_\alg\mid \alg\in \Rg\}$ of root vectors 
such that
$ (E_\alg,E_{-\alg})=1$ for each $\alg\in\Rg_+$,
and put 
$\ol \n_\pm=\+_{\alg\in\Rg_{+}}\ol\g_{\pm \alg}
=\+_{\alg\in\Rg_+}\C E_{\pm \alg}$.
Let $\{H_a\}_{a=1}^{m-1}$ be an orthonormal basis of $\ol\h$.

The classical $r$-matrix of $\ol\g$ is defined by
\begin{equation}
  \label{crm}
  \ol r=\frac{1}{2}\sum_{a=1}^{m-1}H_a \* H_a
+\sum_{\alg\in \Rg_+} E_\alg \* E_{-\alg}\in\ol\g\*\ol\g.
\end{equation}
Choose a set $\{\alg_0,\dots,\alg_{m-1}\}$ of simple roots of $\g$,
then the Weyl group $W_\g$  of $\g$
is generated by the simple reflections $\sigma_a$ 
corresponding to $\alg_a$ $(a=1,\dots,m-1)$.
We put $\rho_{\ol\g}=\frac{1}{2}\sum_{\alg\in R_{\ol\g +}}\alg$,
and denote its dual by $\rho\ch_{\ol\g}\in\ol\h$.
Put
$$\begin{array}{l}
  \rho_\g=\rho_{\bar\g}+mc^*_\g,\quad
\rho\ch_\g=\rho\ch_{\bar\g}+m d_\g.
\end{array}$$

We have a triangular decomposition
$\g=\n_+\+\h\+\n_-$
 with $\n_\pm=\ol\n_\pm\+(\ol\g\*\C[t^{\pm 1}]t^{\pm 1})$.
and put $\b_{\pm}=\n_{\pm}\+ \h$.
Define the subalgebra $\g'$ of $\g$, and its subalgebra $\h'$ and 
$\b_\pm'$ by
 \begin{equation}
    \g'=[\g,\g]=\ol\g\*\C[t,t^{-1}]\+\C c_\g,\quad
 \h'=\ol\h\+\C c_\g,\quad \b_\pm'=\n_\pm\+\h'.
\end{equation}
\subsection{The universal classical $r$-matrix}
We introduce the universal classical $r$ matrix of $\g$ and
some of its properties, which are used in the next section.
Let $\{J_k \}_{k=1,\dots,\dim\ol\g}$ be an orthonormal basis of
 $\ol\g$.
\begin{df}
The universal classical $r$-matrix $r$ of $\g$ is defined by
 \begin{equation} \label{eq;ucrm}
    r=\ol r+\frac{1}{2}(c_\g\* d_\g+ d_\g\* c_\g)+
    \sum_{n\geq 1}\sum_{k=1}^{{\rm dim}\ol\g}
    (J_k\*t^n)\* (J_k\* t^{-n}).
     \end{equation}
\end{df}
For each $X\in\g$, we put 
\begin{equation}
  \ba (X)=[r,X\*1+1\* X]\quad (X\in\g).
\end{equation}
Although the $r$-matrix itself  
contains an infinite sum, $\ba(X)$ always belongs to $\g\*\g$, 
and thus
we have the map $\ba:\g\to\g\*\g$ called the Lie bi-algebra 
structure of
$\g$.
We let $\trans$ be the transposition on $\g\*\g$
and $\proj$ be the natural projection from $\g\*\g$ to $U(\g)$:
\begin{equation}
  \trans(X\* Y)=Y\* X ,\quad\proj(X\* Y)=XY \quad(X,Y\in\g).
\end{equation}
The following two lemmas play important roles in the next section and
can be shown
by direct calculations: 
\begin{lem}For $X\in\g$,
 we have
  \begin{eqnarray}
    \trans\ba(X)&=&-\ba(X),\label{eq;trans}\\
    \proj \ba(X)&=&[\rho\ch_\g \,,\, X].
\label{eq;proj}
  \end{eqnarray}
\end{lem}
For a (possibly infinite) sum $x=\sum_k X_k\* Y_k$ $(X_k,Y_k\in\g)$, 
we set
$x_{1,2}=x\* 1$, $x_{2,3}=1\* x$ and $x_{1,3}=\sum_k X_k\*1\* Y_k$.
\begin{lem}\label{lem;cyb1}
The classical Yang-Baxter equation holds$:$
  \begin{equation}\label{eq;cyb1}
    [r_{1,2},r_{2,3}]+ [r_{1,3},r_{2,3}]+ [r_{1,2},r_{1,3}]=0.\Box
  \end{equation}
\end{lem}
\subsection{Representations of $\g$}
Recall the category $\Og$ of left $\g$-modules: its object set 
(which will be denoted also by $\Og$) 
consists of left $\g$-modules which are
finitely generated over $\g$, $\n_+$-locally finite,
and
$\h$-diagonalizable with finite dimensional weight spaces. 
We also define
the lowest version $\O^{\dagger}$ of $\Og$:
 its objects are
finitely generated over $\g$, $\n_-$-locally finite,
and
$\h$-diagonalizable with finite dimensional weight spaces.
For a weight $\lm\in\h^*$, let $M(\lm)$ denote the
highest-weight Verma module with highest-weight $\lm$, and
$M^{\dagger}(\lm)$ the
lowest-weight  Verma module with lowest-weight $-\lm$.
Their irreducible quotients are denoted by $L(\lm)$ and
 $L^{\dagger}(\lm)$
respectively.
Apparently $M(\lm),L(\lm)\in\Og$ and
$M^{\dagger}(\lm),L^{\dagger}(\lm)\in\Og^{\dagger}$.

For $\l\in\C$, a $\g$-module is said to be  of level $\l$ if
 the center
$c_\g$
acts as the scalar $\l$.
In the following we fix a complex number $\l$.
Let $\Og(\l)$ (resp. $\Og^{\dagger}(\l)$) be the full subcategory
of $\Og$ (resp. $\Og^{\dagger}$) whose object set consists of level 
$\l$ (resp. $-\l$) objects in $\Og$ (resp. $\Og^{\dagger}$).

Put $\h^*(\l)=\{\lm\in\h^*\mid \lm(c_\g)=\l\}$ and
let $\Pg(\l)$ (resp. $\Pg^+(\l)$) be the subset of $\h^*(\l)$ 
consisting of
integral (resp. dominant integral) weights. 
Note that $L(\lm)$ and $L^{\dagger}(\lm)$ become integrable for
$\lm\in\Pg^+(\l)$
and that $\Pg(\l)$ is empty unless $\l$ in an integer.

We also introduce the evaluation module, which do not belong to
 $\O$ nor
$\O^{\dagger}$.
Let $\Vbox=\C^{\m}$ be the vector representation of $\ol\g$
with natural basis $\{u_a\}_{a=1}^m$.
It is convenient for later use to introduce notations 
$\ol e_a\in\bar\h^*$
$(a=1,\dots,m)$ 
denoting the weight of $u_a$. Note that
$\sum_{a=1} ^m\ol e_a=0$.
On the space $\C[z,z^{-1}]\*\Vbox$, 
we define a $\g$-module structure of level $0$ from 
the correspondence
$$ X\*f(t)\mapsto f(z)\* X,\
 d_\g\mapsto z\frac{\partial}{\partial z}\* id.$$
We identify the space $\C[\ol P]\*\VN$ with 
$\bigotimes_{i=1}^N\left(\C[z_i,z_i^{-1}]\* \Vbox\right)$, 
on which $\g$
acts diagonally.

For a $\g$ (resp. $\ol \g$)-module $V$ and a weight $\lm\in\h^*$(resp.
$\ol\h^*$)
let $V_\lm$ denote the weight space of $V$ of weight $\lm$.
For $\lm\in\h^*(\l)$, the image of $\lm$
under the projection $\h^*\to\ol\h^*$ is denoted by $\ol\lm$ and
called the classical part of $\lm$.
\section{Construction of $\H$-modules}
Throughout this paper
we use the notation
\begin{equation}
  M/{{\frak a}}=M/{{\frak a}}M,
\end{equation}
for a Lie algebra ${{\frak a}}$ and an $ {{\frak a}}$-module $M$.
For $A\in\Og(\l)$ and $B\in\Og^{\dagger}(\l)$,
we put
\begin{equation}  \label{eq;bifunctor}
 \F(A,B)=\left.\left(A\*\C[\ol P]\*\VN \* B\right)\right/
  \g',
 \end{equation}
where $\g'$ acts diagonally on the numerator.
We will  define an action of the 
degenerate double affine Hecke algebra $\H$
on the space $\F(A,B)$
using  the expression $\H=\C[W]\* S[\t']$
(see Proposition \ref{eq;another}).
\subsection{KZ connections}
We first define an action of $\bigW=\smallW\ltimes \smallP$
as follows;
define the action of $\ol W$ on  $A\*\C[\ol P]\*\VN\* B$
by $\tau=\text{id}_A\*\tau_1\*\tau_2\*\text{id}_B$, where
$\tau_1$ denote the action on  $\C[\ol P]$ and 
\begin{equation}
    \tau_2(w)(v_1\*\cdots\*v_N)=v_{w(1)}\*\cdots\* v_{w(N)},
\end{equation}
for $w\in\ol W,$ and $v_1,\dots,v_N\in \Vbox.$
Let the elements of $\C[\ol P]$ act on $A\*\C[\ol P]\*\VN\* B$
by usual multiplication.
Apparently this action 
defines an action of $\bigW$ on $\F(A,B)$.

To define an action of $S[\t']$,
we introduce the KZ-connections in terms of the 
universal classical $r$-matrix of $\g$.
Since they do not preserve  $\F (A,B)$ or $A\*\C[\ol P]\*\VN\* B$,
we need to consider localizations of these spaces once.
Consider the manifold
\begin{equation}
  X:=(\C^*)^N\setminus\{~ (z_1,\dots,z_N)\in (\C^*)^N\mid
z_i=z_j\hbox{ for some }i\neq j~\}.
\end{equation}
Let $\OX$ be the sheaf of holomorphic functions on $X$
and let 
$\OX(U)$ is the ring of holomorphic functions on 
an open submanifold $U$ of $X$.
We regard
$\C[\ol P]\*\Vbox^{\* N}$
as a subspace of $\OX(X)\*\VN$.
The diagonal action of $\g$ on
$A\* \C[\ol P]\*\VN\* B$
extends naturally to the one on $A\* \OX(X)\*\VN\* B$.

 Fix a pair of $\g$-modules $A,B$ such that
$A\in\Og(\l)$, $B\in\Og^{\dagger}(\l)$.
For $i\in\{0,\dots,N+1\}$,
let 
$\act_i : \g\to\g^{\* (N+2)}$ 
be the embedding to the $(i+1)$-th component.
Then $\act_i$ induces the map
$\g\to\End_\C( A\* \C[\ol P]\*\Vbox\* B)$ 
and is extended to the map
$\g\to\End_\C(A\*\OX(X)\*\VN\* B)$, which are  denoted by
 the same symbol
$\act_i$.
We set $\ij=\act_i\*\act_j$
and put
\begin{equation}
    r_{(0\,i)}=\act_{0\,i}(r),\  
    r_{(i\,N+1)}=\act_{i\,N+1}(r)\quad  (1\leq i\leq N),
\end{equation}
which are well-defined operators on $A\*\OX(X)\*\VN\* B$.
On the other hand,
 for $i,j\in\{1,\dots,N\}$ with $i<j$, 
we have formally
\begin{equation}
\ij(r)=
\ij(\ol r)+\sum_{n\geq1}(z_i/z_j)^n\ij(\Omega),
\end{equation}
where $\Omega=\ol r+\trans(\ol r)$.
These operators converge simultaneously only on
the region 
 $X_0=\{~(z_1,\dots,z_N)\in X\mid |z_1|<|z_2|<\cdots<|z_N|~\}$ of $X$
to the elements
\begin{equation}\label{eq;rij}
  r_{(ij)}=
  \ij(\ol r)+\frac{z_i/z_j}{1-z_i/z_j}\ij(\Omega).
\end{equation}
of $\End_\C(A\*\OX(X_0)\*\VN\* B)$. 
We extend  $r_{(i,j)}$ holomorphically
to an element in $\End_\C(A\*\OX(X)\*\VN\* B)$.
Then Lemma \ref{lem;cyb1} implies the following:
\begin{lem} \label{lem;cyb2}
The family of operators $\{r_{(ij)}\}_{0\leq i<j\leq N+1}$ satisfies
the classical Yang Baxter equations{\rm :}
\begin{equation}\label{eq;cyb2}
  [r_{(ij)},r_{(jk)}]+ [r_{(ik)},r_{(jk)}]+ [r_{(ij)},r_{(ik)}]
=0
\ \ (0\leq i<j<k\leq N+1).
\end{equation}
\end{lem}
\begin{df}
We define the KZ-connection 
$\nabla_i$ $(i=1,\dots,N)$ as an element of 
$\End (A\*\OX(X)\*\VN\* B)$ 
given by 
\begin{equation}
  \label{eq;KZ}
  \nabla_i=\sum_{0\leq j<i}r_{(ji)}-\sum_{i<j\leq N+1}r_{(ij)}
  +\act_i\left(\rho\ch_\g\right)
\end{equation}
\end{df}
The following proposition follows from direct calculations 
using Lemma \ref{lem;cyb2}.
\begin{pr}\label{pr;KZproperties}$ $
The followings  hold in $\End (A\*\OX(X)\*\VN\* B):$
  \begin{eqnarray}
    &[\nabla_i,\nabla_j]=0 &(1\leq i,j\leq N), \label{eq;commute}\\
    &\tau(w)\, \nabla_i\, \tau(w^{-1}) = \nabla_{w(i)} &(w\in\ol W),
\label{eq;wnabla}\\
    &\left[ \nabla_i , f\right]        = (\l+\m) \partial_i f
&(f\in\OX(X)).
\label{eq;nablaf}
  \end{eqnarray}
\end{pr}
\begin{pr}
\label{pr;KZpreserve}
Put $\act=\sum_{i=0}^{N+1}\act_i.$
Then
we have 
\begin{equation}\label{eq;KZX}
  [\KZ_i,\act(X)]=(\act\*\act_i)(\ba(X))
\end{equation}
for each $i=1,\dots,N,$ and $X\in\g$,
where $\act_i\*\act_i(\ba(X))$ means $\act_i\proj(\ba(X))$.
In particular, the KZ connections $\nabla_i$ 
preserve 
the subspace $\g'(A\*\OX(X)\*\VN\* B):$
  \begin{equation}
    \nabla_i \left(\g'(A\*\OX(X)\*\VN\* B)\right)
    \subset \g'(A\*\OX(X)\*\VN\* B).
  \end{equation}
\end{pr}
\noindent
{\it Proof.}
For $X\in\g$, it is obvious that
$[r_{(ij)},(\act_i+\act_j)(X)]=\ij(\ba(X))$,
and thus we have
$$\begin{array}{l}
[\nabla_i,\sum_{k=0}^{N+1}\act_k(X)]
=\sum_{0\leq j<i}\ji(\ba(X))-\sum_{i<j\leq N+1}
\ij(\ba(X)) +\act_i [\rho\ch_\g , X]\\
=\sum_{0\leq j\leq N+1,j\neq i}\ji(\ba(X))
+ \act_i\varpi(\ba(X))
=\sum_{0\leq j\leq N+1}\ji(\ba(X)).
\end{array}$$
The second equality
follows from (\ref{eq;trans}) and (\ref{eq;proj}). 
\qed
\subsection{The Cherednik-Dunkl operator}
We generalize the Cherednik-Dunkl operators
by combining with the KZ connections
to the operators on $A\*\OX(X)\*\VN\* B$, which turn out to act
on $\F(A, B)$.

For each $\kappa\in\C$ and $\xi\in\bar\t$, 
the Cherednik-Dunkl operator $D_\xi=D_\xi(\kappa)\in\End_\C(\OX(X))$ 
is given by 
  \begin{equation}
    D_\xi=
\kappa\,\partial_\xi+\sum_{\al\in \ol R_+}\al(\xi)
\frac{1-\tau_1(s_\al)}{1-e^{-\al}}
    -\ol\rho(\xi),
  \end{equation}
where $\ol\rho=\frac{1}{2}\sum_{\al\in \ol R_+}\al\in\bar\t^*$.
Note that $D_\xi$ preserves $\C[\ol P]$. 
\begin{pr}[Cherednik]\label{pr;Dunkl}
  The correspondence
\begin{equation}
  \begin{array}{lll}
    c  &\mapsto \kappa,      & \\
    \xi&\mapsto D_\xi\ \     &(\xi\in\bar\t),\\
    w  &\mapsto \tau_1(w)\ \ &(w\in\ol W),\\
    f  &\mapsto f\times \    &(f\in\OX(X))
  \end{array}
\end{equation}
defines an action of $\H$ on $\OX(X)$ and on $\C[\ol P]$.
\end{pr}
Let us generalize the Cherednik-Dunkl operators to our case.
For $\xi\in\bar\t$, we set 
\begin{equation}
  \nabla_\xi=\sum_{i=1}^N \e_i(\xi)\nabla_i
\end{equation}
and define the operator $\D_\xi=\D_\xi(\k)\in\End(A\*\OX(X)\*\VN\*B)$ 
by
  \begin{equation} \label{eq;genDunkl}
    \D_\xi:=\nabla_\xi+
    \sum_{\al\in \ol{R}_+}\al(\xi)\frac{1-\tau(s_\al)}{1-e^{-\al}}
    -{\ol\rho}(\xi)-\frac{m-1}{2m}\partial_\xi\hbox{log}G,
  \end{equation}
where 
 $  G=\Pi_{\al\in \ol R}(1-e^{-\al})\in\C[\ol P].$

\begin{th} \label{th;Haction}
  $(${{\rm i}}$)$  The correspondence
  \begin{equation}\label{eq;Haction}
      \begin{array}{lll}
      c&\mapsto \l+m,&\\  
    \xi&\mapsto \D_\xi\   &(\xi\in\bar\t),\\
    w  &\mapsto \tau(w)  &(w\in\ol W),\\
    f  &\mapsto f\times  &(f\in\OX(X))
      \end{array}
  \end{equation}
defines an action of $\H$ on 
$A\*\OX(X)\*\VN\* B$. 

$(${{\rm ii}}$)$ 
The subspaces
$A \* \C[\ol P] \* \VN\* B$ and 
$\g'(A\*\C[\ol P]\*\VN\* B  )$ are preserved by the above action.
Therefore the correspondence $(\ref{eq;Haction})$ induces 
an action of $\H$
on the space $\F(A,B)$.
\end{th}
\noindent
{\it Proof.}
(i) By Proposition \ref{pr;KZproperties}, the relations between the
operators 
$\nabla_\xi$, $w$ and $f$ 
are same with the ones between $(\l+m)\partial_\xi$,
$w$ and $f$. 
Therefore the statement
is shown exactly as Proposition \ref{pr;Dunkl}.

(ii)
The fact that $\D_\xi$ preserves $A \* \C[\ol P] \* \VN\* B$
follows from direct calculations.
The rest is straightforward from the definition of the action
and Proposition \ref{pr;KZpreserve}.
\qed

\medskip
In the rest of this paper, we put $\k=\l+m$.
It is easily checked that our construction is
 functorial in the following
sense:
for $\g$-homomorphisms $f_1:A\to A'$ and $f_2:B\to B'$ between
$\g$-modules,
the induced maps
\begin{equation}
  f_{1*}:\F(A,B)\to\F(A',B),\quad f_{2*}:\F(A,B)\to\F(A,B')
\end{equation}
are both $\H$-homomorphisms. Therefore
we have the bifunctor from $\Og(\l)\times\Og^{\dagger}(\l)$
to the category of $\H(\k)$-modules.
\section{Structure of $\H$-modules}
Our next aim is to study $\H$-module structures of
$\F(A,B)$ for $A\in\Og(\l)$ and $B\in\Og^{\dagger}(\l)$.
For each $\mu\in \h^*(\l)$, we can define the right exact functor
 \begin{equation}
    \F_\mu(\cdot)=\F(M(\mu),\cdot\,)
  \end{equation}
from $\Og^{\dagger}(\l)$ to the category of $\H$-modules.
We put 
$$
  \label{eq;FMM}
         \M(\pair)=       \F_\mu(M^{\dagger}(\lm)),\quad
         \V(\pair)=       \F_\mu(L^{\dagger}(\lm)),
$$
 for each $\pair\in\h^*(\l)$.
When $\lm,\mu\in\Pg^+(\l)$, the space $\V(\pair)$
is related to  the dual space of conformal blocks (see Proposition
\ref{pr;propagation} and \cite{TUY}).
\subsection{Fundamental properties}
In this subsection,
 we study some properties of the bifunctor $\F$, which
will be used later.
Let $v(\mu)$ and $v^{\dagger}(\lm)$ denote the highest (and lowest) 
weight vector
of $M(\mu)$ and $M^{\dagger}(\lm)$ respectively.
\begin{lem}\label{lem;restrict}
For $A\in\Og(\l),B\in\Og^{\dagger}(\l),\lm\in \h^*(\l)$, 
and $\mu\in\h^{*}(\l)$,
the natural inclusions $\C v(\mu)\to M(\mu)$ and
$\C v^{\dagger}(\lm)\to M^{\dagger}(\lm)$ induce isomorphisms
\begin{eqnarray}
   \F(M(\mu),B)  &\cong&
   \C v(\mu)\*\C[\ol P]\*\VN \* B/
    \b_+', \\  
   \F(A,M^{\dagger}(\lm))&\cong&
    A\*\C [\ol P]\*\VN \* \C v^{\dagger}(\lm)/
    \b_-'.  
    \nonumber
\end{eqnarray}
\end{lem}
\noindent
{\it Proof.}
Follows from the following well-known $\g$-isomorphism 
for any $\g$-module $V$:
\begin{equation}\label{eq;Mackey}
  M(\mu)\* V\cong U(\g)\*_{U(\b_+)} (\C v(\mu)\* V)
\end{equation}
and its lowest version, where 
on the left-hand-side of $(\ref{eq;Mackey})$, $\g$ acts 
diagonally and
on the right-hand-side by left multiplication.
\qed

\medskip
Again by  $(\ref{eq;Mackey})$,
there is an isomorphism of vector spaces
$$\M(\pair)\cong \C v(\mu)\*(\C[\ol P]\*\VN)_{\lm-\mu}\*\C
v^{\dagger}(\lm).$$ 
Notice that the space on the right-hand-side 
is an $\bigH$-submodule of 
$M(\mu)\*(\C[\ol P]\*\VN)\* M^{\dagger}(\lm)$
with respect to the action  $(\ref{eq;Haction})$.
Hence we have
\begin{pr}\label{pr;Miso}
  For $\lm,\mu\in \h^*(\l)$, 
we have
  \begin{equation}
\M(\pair)\cong \C v(\mu)\*(\C[\ol P]\*\VN)_{\lm-\mu}\*
\C v^{\dagger}(\lm)
  \end{equation}
as $\bigH$-modules.
\end{pr}

\begin{pr}\label{pr;propagation}
For $\lm,\mu\in \Pg^+(\l)$, the natural projection
$M(\mu)\to L(\mu)$ induces
an $\H$-isomorphism
\begin{equation}
 \F_\mu(L^{\dagger}(\lm))\cong \F (L(\mu),L^{\dagger}(\lm)).
\end{equation}
\end{pr}
{\it Proof.}
It is obvious that the induced map 
$\F(M(\mu),L^{\dagger}(\lm)) \to \F(L(\mu),L^{\dagger}(\lm))$
is surjective.
To show the injectivity, recall that the maximal submodule of 
$M(\mu)$ is 
generated by the singular vectors 
$v_a\in M(\mu)_{\sigma_a\circ\mu}$ $(a=0,\dots,m-1)$ 
satisfying $X v_a=0 $
for any
$X\in\n_+$, where 
$\sigma_a\circ\mu=\sigma_a(\mu+\rho_\g)-\rho_\g$.
By elementary calculations, it 
can be checked that, for any number $n\geq0$,
there exists a non-zero constant $c_n$ such that 
$E_{\alg_a}^n E_{-\alg_a}^n v_a=c_n v_a$, 
where $E_\alg$ denotes the root 
vector corresponding to $\alg\in\R_g$.
Put $V=\C[P]\*\VN\* L^{\dagger}(\lm)$. 
It is enough to show that
$U(\g)v_a\* V\subset\g(M(\mu)\* V)$
for each $a=0,\dots,m$, which is reduced to showing 
that
\begin{equation}\label{eq;reduced}
   v_a\* u\in\g(M(\mu)\* V)
\end{equation}
for any $u\in V$.
Since $V$ is integrable and thus $\n_+$-locally finite,
there exists a number $n>0$ such that $E_{a_i}^nu=0$,
and we have
$$c_n  v_a\* u=v_a\* E_{\alg_a}^n E_{-\alg_a}^n u\\ \equiv
 E_{-\alg_a}^n E_{\alg_a}^n  v_i\* u\equiv 0 
\quad\text{mod}\;\g(M(\mu)\* V),
$$
that implies $(\ref{eq;reduced})$.
\qed
\subsection{Isomorphisms to induced modules}\label{ss;induced}
We shall describe the $\H$-module structure
of $\M(\pair)$, which will be shown to be isomorphic to an
 induced module 
$\Y_\beta(\zeta)=\H(\k)\otimes_{\ol\H_\beta}\C \one_{\zeta}$ 
(see \S \ref{ss;standard})
associated to an appropriate parameter $(\beta, \zeta)$ determined by
$\lm,\mu
\in\h^*(\l)$.

Recall that $\bar e_a\in \bar\h^*$ $(a=1,\dots,m)$
denotes the weight of the standard basis $u_a$ of $\Vbox$. 
Let us associate 
an ordered partition $\beta$ of $N$ and an element $\zeta\in (\t')^*$
to each pair $\lm,\mu\in\h^*$.
Assume that $\lm, \mu\in \h^*(\l)$ 
satisfy $\ol\lm-\ol\mu\in\wt(\VN)$.
Then $\ol\lm-\ol\mu=\sum_{a=1}^m \beta_a\ol e_a$ for some
$\beta_a\in\Z_{\geq0}$ $(a=1,\dots,m)$ satisfying 
$\sum_{a=1}^m \beta_a=N$.
Thus we define an ordered partition $\beta_{\pair}$ of $N$ by
\begin{equation}\label{eq;beta}
  \beta_{\pair}
=(\beta_1,\dots,\beta_m).
\end{equation}
Here, we allowed
 appearance of zeros as elements of an ordered partitions.  

We keep the assumption that
$\bar\lm-\bar\mu \in \text{wt}(\VN)$.
Put $\mu_a=(\mu,\bar e_a)$ $(a=1,\dots,m)$
and consider the set 
\begin{equation}
  Y_{\pair}=Y_{\clpair}:=\{~(a,\y)\in \Z\times \C \mid 
  a=1,\dots,m ;\; \y=\mu_a+1,\mu_a+2,\dots,
\mu_a+\beta_a~\},
\end{equation}
which we represent just by $Y$ below.
A bijection $T :Y\to\{1,\dots,N\}$ is called
a tableau on $Y$. 
We let $\T(\clpair)$ denote
the set of tableaux on $Y$.
Define $T_0\in\T(\pair)$ by
\begin{equation}\label{eq;basic1}
T_0(a,\y)=\sum_{k=1}^{a-1}\beta_k+\y-\mu_a.
\end{equation}
Now, define $\ol\zeta_{\clpair}\in \bar\t^*$
and $\zeta_{\pair}\in (\t')^*$
by
\begin{eqnarray}\label{eq;zeta}
   \ol\zeta_{\clpair}
   &=&\sum_{(a,\y)\in Y}(\y-a)\e_{T_0(a,\y)}
   +\frac{1}{2m}(m^2-N)\sum_{k=1}^N \e_k\\
\label{eq;doublezeta}
  \zeta_{\pair}&=& \ol\zeta_{\clpair}+\k c^*
\end{eqnarray}
For $T\in\T(\clpair)$, let
$w_T$ denote 
the inverse image of $T\in\T(\clpair)$
under the isomorphism
$$\begin{array}{rcl}
  \ol W&\to& \T(\clpair)\\
  w&\mapsto&           w T_0:(a,\y)\mapsto w(T_0(a,\y)).
\end{array}$$
Note that
\begin{equation}\label{eq;wzeta}
   w_T(\ol\zeta_{\clpair})
   =\sum_{(a,\y)\in Y}(\y-a)\e_{T(a,\y)}
   +\frac{1}{2m}(m^2-N)\sum_{k=1}^N \e_k.
\end{equation}
The following lemma can be checked easily.
\begin{lem}
  If $\bar\lm-\bar\mu\in\text{{{\rm wt}}}(\VN)$, then
$\zeta_{\pair}$ $($resp. $ \ol\zeta_{\clpair})$
belongs to $(\t')^*_{\beta_{\pair}}$ $($resp. 
$\bar\t^*_{\beta_{\pair}})$
$($see $(\ref{eq;tbeta})$ for the definition of $(\t')^*_\beta$ $)$.
\end{lem}
Let
$u(\mu)$ and $u^{\dagger}(\lm)$ be the highest and
 lowest weight vector in 
$M(\mu)$ and $M^{\dagger}(\lm)$ respectively, and put
\begin{equation}
  u_{\pair}=v(\mu)\*\overbrace{u_1\*\cdots\* u_1}^{\beta_1}
\*\overbrace{u_2\*\cdots\* u_2}^{\beta_2}\*\cdots\*
\overbrace{u_m\*\cdots\* u_m}^{\beta_m}\* v^{\dagger}(\lm).
\end{equation}
By direct calculations, we have the following lemma.
\begin{lem}\label{lem;weight}
  For $\lm,\mu\in \h^*(\l)$ such that $\bar\lm-\bar\mu\in 
\text{{{\rm wt}}}(\VN)$,
  the vector $u_{\pair}$ is
 a weight vector with the weight $\zeta_{\pair}:$
  \begin{equation} 
    \D_\xi u_{\pair}=\zeta_{\pair}(\xi)
    u_{\pair}\quad \text{for }\xi\in\t.
  \end{equation}
\end{lem}
By Proposition \ref{pr;Miso} and Lemma \ref{lem;weight},
we get the following conclusion:
\begin{pr}\label{pr;standard}
For any $\lm,\mu\in\h^*(\l)$, we have a
unique $\H$-module isomorphism 
\begin{equation}
    \M(\pair)\cong
\left\{
\begin{array}{ll}
\Y_{\beta_{\pair}}(\zeta_{\pair})\quad &\text{ if }\bar\lm-\bar\mu\ 
\in \text{{{\rm wt}}}(\Vbox^{\*N}),\\
0& otherwise,
\end{array}
\right.
\end{equation}
which sends 
$u_{\pair}$ to $\one_{\zeta_{\pair}}$.
\end{pr}
%
\subsection{Finite analogue of the functor $\F$}\label{ss;finite}
It is possible to construct
the ``finite version'' of our functor $\F$,
which is easier to study and has some suggestive properties.
Consider the category $\ol\O$ (resp. $\ol\O^*$)
of highest (resp. lowest) weight $\ol\g$-modules. 
For $\ol\g$-modules $A\in\ol\O$ and $B\in\ol\O^*$, we put
\begin{equation}
  \ol\F(A,B):=\left.\left( A\*\VN\* B\right)\right/\ol\g.
\end{equation}
and consider the operator 
\begin{equation}
  \ol\D_i :=\sum_{0\leq j< i}\left( \ji(\ol{r})+\frac{1}{2m}\right)-
\sum_{i<j\leq N+1}\left(\ij(\ol{r})+\frac{1}{2m}\right)+\act_i(\rhog).
\end{equation}
It can be shown that $\ol\D_i$ acts on $\ol\F(A,B)$ and
the correspondence
\begin{equation}\begin{array}{lll}
\e_i\ch &\mapsto \ol\D_i  \quad &i=1,\dots,N,\\
w       &\mapsto \tau_2(w)\quad &w\in\ol W 
\end{array}\end{equation}
defines an action of $\ol\H$ on $\ol\F(A,B)$.
For each $\mu\in\ol\h^*$, we
 define a functor from the category $\ol\O^*$ to
the category of finite dimensional $\ol\H$-modules by
\begin{equation}
  \ol\F_{\bar\mu}(\cdot)=\ol\F(\ol M(\bar\mu),\cdot),
\end{equation}
and put  
$\ol\F_{\bar\mu}(\ol M^{\dagger}(\bar\lm))=\ol \M(\clpair),\quad
 \ol\F_{\bar\mu}(\ol L^{\dagger}(\bar\lm))=\ol \V(\clpair),
$
where $\ol M(\bar\lm)$ denotes the highest weight left Verma module
of $\ol\g$ with highest weight $\bar\lm\in\ol\h^*$ etc.
As an analogue of Proposition \ref{pr;propagation}
we have $\ol \V(\clpair)\cong\ol\F(\ol L(\bar\mu),\ol
L^{\dagger}(\bar\lm))$
if $\ol\lm,\ol\mu\in\ol\h^*$ are both dominant integral.
An analogue of Proposition
\ref{pr;standard} also holds:
\begin{pr}
For any $\bar\lm,\bar\mu\in\ol\h^*$,
we have the following $\ol\H$-module isomorphism:
\begin{equation}
   \ol\M(\clpair)\cong
\left\{
\begin{array}{ll}
\ol\Y_{\beta_{\clpair}}(\ol\zeta_{\clpair})
\quad &\text{ if }\bar\lm-\bar\mu\
\in \text{{{\rm wt}}}(\Vbox^{\*N}),\\
0& otherwise,
\end{array}
\right.
\end{equation}
where $\beta_{\clpair}$ and $\ol\zeta_{\clpair}$
 are defined by the similar formulas as $(\ref{eq;beta})$
and $(\ref{eq;zeta})$ respectively.
\end{pr}

Under finite situations, standard arguments using the infinitesimal
characters 
 deduce the following:
\begin{pr}\label{pr;exact}
 Suppose that $\bar\mu\in\ol\h^*$ satisfies 
$$(\bar\mu+\rho_{\ol\g},\alg)\notin\{-1,-2,\dots\}\quad
\text{for all } \alg\in
R_{\ol\g+},$$
then
the functor $\ol\F_{\bar\mu}$
is exact.
\end{pr}
Now, let us suppose that $\ol\lm,\ol\mu$ are
both dominant integral.
Then $\ol L(\ol\lm)$ is finite dimensional and
the dimension of $\ol\F_{\ol\mu}(\ol L(\ol\lm))$ is
calculated 
 using Shur-Weyl reciprocity law and
turns out to equal the number of standard tableaux
on the skew Young diagram $Y_{\clpair}$:
A tableau $T\in\T(\clpair)$ on  $Y_{\clpair}$
is called a standard tableau if $T$ satisfies the following two
conditions:
$$\begin{array}{l}
T(a,\y)<T(a,\y+1) 
\quad\text{if }(a,\y),(a,\y+1)\in Y_{\clpair},\\
T(a,\y)<T(a+1,\y)
\quad\text{if }(a,\y),(a+1,\y)\in Y_{\clpair}.
\end{array}$$
Let $\T_s(\clpair)$
denote the set of standard tableaux on $Y_{\clpair}$ 
(then $\dim\V(\clpair)=\# \T_s(\clpair)$),
Note that $w_T\in \bar W^{\beta_\pair}$ for $T\in\T_s(\pair)$.
The following proposition can be proved by using 
Proposition \ref{int:lead} and Proposition \ref{pr;w(zeta)=zeta}.
\begin{pr}\label{pr;tame}
  Let $\bar\lm,\bar\mu\in\ol\h^*$ be both dominant integral.
Then
$\ol\V(\clpair)$ is the unique irreducible quotient of
$\ol\M(\clpair)$. Moreover, its basis is given by
$\{\varphi_{w_T}{\one}_{\zeta_{\clpair}}\}_{T\in\T_s(\clpair)}$.
In particular $S[\,\ol\t\,]$ acts semisimply
on $\ol\V(\clpair)$.
\end{pr}
\subsection{Dominant integral case}\label{ss;dominant_integral_case}
As a consequence of Proposition \ref{pr;standard},
we have that the $\H$-module $\M(\pair)$ is irreducible if
 $\zeta_{\pair}$ is $\beta_{\pair}$-generic 
(Definition \ref{df;generic}).

Let us consider the case where 
$\lm$ and $\mu$ are both dominant integral.
The following statement is a corollary of Corollary
\ref{cor;sufficient}: 
\begin{cor}\label{cor;V}
  If $\mu\in \Pg^+(\l)$,
then $\M(\pair)$ has a unique irreducible quotient.
\end{cor}
{\it Proof.}
It can be checked that the condition $\mu\in \Pg^+(\l)$
implies the condition $(\ref{cond;dominant})$ 
in Corollary \ref{cor;sufficient} for $\zeta_{\pair}$.
\qed
\begin{con}
\label{con;tame}
  Let  $\lm,\mu\in\h^*(\l)$ be both dominant integral.
Then
$\V(\pair)$ is the unique irreducible quotient of
$\M(\pair)$. Moreover $S[\,\t']$ acts semisimply
on $\V(\pair)$.
\end{con}
%

\section{Symmetric part}
In the rest of this paper
we focus on the $\ol W$-invariant part of $\V(\pair)$ for
$\lm,\mu\in\Pg^+(\l)$
and attempt to give a description corresponding 
to Corollary \ref{cor;sympart}. 
%
\subsection{Specialized characters}
As an approach to investigate the $\H$-modules, we
 consider their characters with respect to the operator
\begin{equation}
   L_1:=\frac{1}{\l+m}\left\{
\SumN \D_i-\frac{1}{2m}N(m^2-N)\right\}
\in S[\,\ol\t\,]^{\ol W}.
\end{equation}
To this end, we introduce ``the polynomial part'' of 
$\H$ and $\F_\mu(B)$
(see \S \ref{ss;DDAHA}).
Put
\begin{equation} \F_\mu\positive  (B)=\C v(\mu)\*
\C[\ol P\positive ]\*\VN\*
B/\b_+'.
\end{equation}
Then $\F_\mu\positive (B)$ has 
the $\H\positive $-module structure.
For a vector space $M$ with a semisimple action of  $L_1$
with a finite dimensional eigen space decomposition, we put
\begin{equation}
  \char M = \hbox{Trace}_M q^{L_1}
\end{equation}
and call 
$\char M$ the specialized character of $M$, 
where $q$ is a parameter.
It is easily checked that 
$\char \F_\mu\positive (B)$ is well-defined for any
$B\in\Og^{\dagger}(\l)$.
Put $\M\positive (\pair)= \F_\mu\positive (M^{\dagger}(\lm))$ and 
$\V\positive (\pair)= \F_\mu\positive (L^{\dagger}(\lm))$.

The same argument as Lemma \ref{pr;Miso}
shows that 
 $$\M\positive (\pair)\cong
\C v(\mu)\*(\C[\ol P\positive ]\*\VN)_{\lm-\mu}\* 
\C v^{\dagger}(\lm),$$
as an $\H$-module.
On the right-hand-side, we have
\begin{equation}\label{eq;L1}
L_1\equiv  \triangle_{\bar\lm}-\triangle_{\bar\mu}+\sum_{i=1}^N
\partial_i,
\end{equation}
where $\triangle_{\bar\lm}=\frac{1}{2(\l+m)}
\left(
(\bar\lm,\bar\lm)+2(\ol\rho_\g,\,
\bar\lm)\right)$.
>From (\ref{eq;L1}), the specialized characters of 
$\M\positive (\pair)$ and 
$\M\positive (\pair)^{\ol W}$
are calculated as
$$\begin{array}{l}
  \char \M\positive (\pair)=
  \left(
    \begin{array}{c}
N\\ \beta_{\pair}
    \end{array}
    \right)
  \frac{q^{\triangle_{\bar\lm}-\triangle_{\bar\mu}}}{(q)_N},\\
  \char \M\positive (\pair)^{\ol W}=
\left[
\begin{array}{c}
N\\
\beta_{\pair}
\end{array}
\right]
 \frac{q^{ \triangle_{\bar\lm}-\triangle_{\bar\mu}}}{(q)_N},
\end{array}$$
respectively, 
where 
$(q)_N=(1-q)(1-q^2)\cdots (1-q^N)$
and 
$$
 \left(
    \begin{array}{c}
N\\ \beta_{\pair}
    \end{array}
\right)
=\frac{N!}{\beta_1!\beta_2!\cdots\beta_m!},\quad
  \left[
\begin{array}{c}
N\\ \beta_{\pair}
\end{array}
\right]
=\frac{(q)_N}{(q)_{\beta_1} (q)_{\beta_2}\cdots (q)_{\beta_m}}. 
$$

It is natural to attempt to construct a resolution of $\V(\pair)$ 
by induced modules to calculate the character for $\V(\pair)$.
 Recall the BGG resolution for $L^{\dagger}(\lm)$, which is the exact
sequence   
\begin{equation}
0\gets L^{\dagger}(\lm)\gets M^{\dagger}(\lm)\gets C_1\gets 
\cdots\gets
C_i\gets\cdots\,
  \label{eq;BGG}
\end{equation}
of $\g$-modules
with 
$$C_i=\+_{
w\in \Wg^{(i)}
}
M^{\dagger}(w\circ \lm),$$ 
where
$\Wg^{(i)}$ is the subset of the Weyl group of $\g$ 
consisting of elements of length $i$.
Sending by the functor $\F_\mu$ (resp. $\F_\mu\positive $),
we have the complex of $\H$ (resp. $\H\positive $)-modules
 respectively:
\begin{eqnarray}{}
  0\gets  \  \V(\pair)     \gets \ \M(\pair)      
\gets  {\cal C}_1  \gets \cdots \gets   {\cal C}_i \gets\cdots ,
  \label{eq;BGG1}\\
    0\gets      \V\positive (\pair)    \gets      \M\positive (\pair)     
\gets      {\cal C}\positive _1      \gets \cdots\gets       {\cal
C}\positive _i  \gets
\cdots ,
  \label{eq;BGG2}
 \end{eqnarray}
where we put
$  \M\positive (\pair)=\F\positive _\mu(M^{\dagger}(\lm))$ and
$  \V\positive (\pair)=\F\positive _\mu(L^{\dagger}(\lm))$.
Note that each ${\cal C}_i$ or ${\cal C}\positive _i$ 
is a direct sum of induced modules:
\begin{equation}
{\cal C}_i        =\+_{w\in \Wg^{(i)}}\M(w\circ\lm,\mu),\quad
{\cal C}\positive _i=\+_{w\in \Wg^{(i)}}\M\positive (w\circ\lm,\mu).
\end{equation}
Therefore if the sequence  
$(\ref{eq;BGG2})$ is exact, we can calculate the character of 
$\V\positive
(\pair)$.
It can be shown the above sequences are exact 
when $\k$ is large enough,
but in general case it is still conjecture:
\begin{con}\label{con;BGGexact}
The sequences $(\ref{eq;BGG1})$ and $(\ref{eq;BGG2})$ are exact.
\end{con}
Note that 
the above conjecture leads to the formula
\begin{equation}
\label{eq;character1}
\hbox{\rm ch}\V\positive (\pair)^{\ol W}
= \sum_{w\in \Wg} (-1)^{l(w)}
q^{(w\circ\lm)(d_\g)} 
\left[
\begin{array}{c}
N\\
\beta_{w\circ\lm,\mu}
\end{array}
\right]
 \frac{q^{ \triangle_{\bar\lm}-\triangle_{\bar\mu}}}{(q)_N},
\end{equation}
\subsection{
Main conjecture}
In the following we fix a pair $(\lm,\mu)\in\Pg^+(\l)\times\Pg^+(\l)$ 
satisfying $\bar\lm-\bar\mu\in\wt(\VN)$.
Recall that
we associated 
a partition $\beta_\pair$ of $N$ such that
$\bar\lm-\bar\mu=\sum_{a=1}^m \beta_a \bar e_a$, and
a diagram
\begin{equation}\label{eq;Young}
  Y_{\pair}=\left\{(a,\y)\in\Z\times \C \mid a=1,\dots,m\, ;\
\y=\mu_a+1,\mu_a+2\dots,\mu_a+\beta_a\right\}.
\end{equation}
where
$\bar\mu_a=(\mu_a, \bar e_i)$.
As in \S \ref{ss;finite},
let $\T(\pair)$ (resp. $\T_s(\pair)$)
denote the set of tableaux (resp standard tableaux) 
on $Y=Y_\pair$.

For $T\in\T(\pair)$ and
a number $n\leq N$, 
let $\beta_T^{(n)}=(\beta_{T\,1}^{(n)},\dots,\beta_{T\, m}^{(n)})$ 
denote the partition of $n$ corresponding to 
the diagram
$Y_T^{(n)}=\{(a,\y)\in Y_{\pair} \mid T(a,\y)\leq n\}
$: namely,
$\beta_{T\, a}^{(n)}=\#\{\y\in\Z+\mu_a\mid (a,\y)\in Y_T^{(n)}\}$.
Put
$$\lm_{T}^{(n)}=
\sum_{a=1}^m\left(\mu_a+\beta_{T\,a}^{(n)}-\frac{n}{m}\right)
\bar e_a+\l
c_\g^*.$$
Then it is checked that $\lm_{T}^{(n)}$ belongs to $\Pg(\l)$. 
A standard tableau $T\in\T_s(\pair)$ is said to be 
$\l$-restricted if
all the weights 
$\lm_T^{(1)}$, $\lm_T^{(2)}$ $,\cdots,$ $\lm_T^{(N)}$
belong to $\Pg^+(\l)$.
We let $\T_s^{(\l)}(\pair)$ denote the set of $\l$-restricted 
standard
tableaux on $Y$.

For $T\in \T(\pair)$, let $w_T\in \bar W^\beta$
be the corresponding element.
Direct calculations imply
 the following lemma:
\begin{lem}
\label{weightcon}
Let $T\in\T_s^{(\l)}(\pair)$.

\noindent
\rm{(i)} For $i\in \{1,\dots,N-1\}$,
let $T(a,\y)=i$ and $T(a',\y')=i+1$.
Then we have
\begin{equation}\label{eq;ineq}
\begin{array}{rllll}
  -(\k-1)\leq &(w_T(\zeta_{\pair}),\al_i)& \leq -2
&\text{ if }&a>a',\\
  &(w_T(\zeta_{\pair}),\al_i)& =-1&\text{ if }&a=a',\\
  1\leq &(w_T(\zeta_{\pair}),\al_i)& \leq\k-2
&\text{ if }&a<a'.
\end{array}
\end{equation}

\noindent
\rm{(ii)}
For any $\al\in\bar R_+$, we have
\begin{equation}
  (w_T\zeta_\pair,\al)\leq\kappa-2.
\end{equation}
\end{lem}
Recall that we associated 
 $\eta_w\in \bar P_-(w)$ for each $w\in \bar W^\beta$
in Definition \ref{df;eta_w}.
\begin{lem}\label{lem;eta}
For each 
$T\in\T_s^{(\l)}(\pair)$, we have
\begin{equation}
\eta_{w_T}=\sum_{i=1}^{N-1} d_{i}(T)(\e\ch_{i+1}+\cdots+\e\ch_N).
\end{equation}
Here 
$$
d_i(T)=\left\{\begin{array}{l}
1 \quad\hbox{ if }a< a'\\
0 \quad\hbox{ if }a\geq a'
\end{array}
\right.,
$$
with $T(a,\y)=i$ and $T(a',\y')=i+1$.
\end{lem}
\begin{rem}
The function $d_i:\T_s^{(\l)}(\pair)\rightarrow 
\N$ can be identified with the so-called
$H$-function in the RSOS model $($see $\cite{Jimbo})$.
\end{rem}
Let $d(T)=\eta_{w_T}(\e_1+\cdots+\e_N)$ and define
\begin{equation}
  \label{eq;F}
   F_{{\pair}}(q)=\sum_{T\in\T_s^{(\l)}({\pair})}q^{d(T)}.
\end{equation}
By the above remark,
the following holds by using standard arguments in the solvable 
lattice model
(see \cite{Jimbo},\cite{Baxter}):
\begin{th}
  \begin{equation}
    F_{{\pair}}(q)= \sum_{w\in \Wg} (-1)^{l(w)}
q^{(w\circ\lm)(d_\g)}
\left[
\begin{array}{c}
N\\
\beta_{w\circ \lm,\mu}
\end{array}
\right].
\end{equation}
\end{th}
Combining with Conjecture \ref{con;BGGexact},
we have a simple formula
\begin{equation} \label{eq;character}
\char \V\positive  (\pair)^{\ol W}
= F_{\pair}(q) 
\cdot
 \frac{q^{ \triangle_{\bar\lm}-\triangle_{\bar\mu}}}{(q)_N}.
\end{equation}
Set $X^{\beta_\pair}_\l=\{ t_\eta w_T\mid
T\in\T_s^{(\l)}(\pair), \ \eta\in \bar P_-(w_T)\}\subseteq 
X^{\beta_\pair}$.
\begin{th}
\label{lower} Let $\lm,\mu\in\Pg^+(\l)$ and
$Q=\frac{1}{N!}\sum_{w\in\bar W}w\in\C[\bar W]$.
 Then
for each $x\in X^{\beta_\pair}_\l$,
the element $Q \varphi_{x}\cdot \one_{\zeta_{{\pair}}}
$ is non-zero.
Moreover the set
$\{ Q \varphi_{x}\cdot \one_{\zeta_{{\pair}}}\}_{
 x\in X^{\beta_\pair}_\l}$ is
linearly independent$:$
$$\V (\mu, \lambda)^{\ol W}\supseteq 
\bigoplus_{x\in X^{\beta_\pair}_\l} \C 
Q \varphi_x
\cdot \one_{\zeta_{{\pair}}} .$$ 
In particular, we have
$$\hbox{{\rm ch}}\V\positive   (\mu, \lambda)^{\bar W} 
\geq  F_{\pair}(q)\cdot
 \frac{q^{ \triangle_{\bar\lm}-\triangle_{\bar\mu}}}{(q)_N}.
$$
\end{th}
{\it Proof.} 
Take any $x\in X^{\beta_\pair}_\l$.
It is easy to see $x\overline{(\zeta_\pair})\in \bar\t^*_{\beta(x)}$.

Direct calculations using Lemma \ref{weightcon} and Lemma 
\ref{lem;eta}
implies that 
\begin{equation}
  \label{eq;con0}
  x\overline{(\zeta_\pair})\in \bar\t^*_{\beta(x)}\text{ is }
\beta(x)\text{-generic},
\end{equation}
and
$
(\zeta_\pair,\al)\neq 0,\pm1\text{ for any }\al\in S(x).
$
In particular,
\begin{equation}
  \label{eq;con2}
  \varphi_x \text{ is invertible on }\one_{\zeta_\pair}.
\end{equation}
By Lemma \ref{lem;Wbetax}, we have
$s_i\varphi_x\one_{\zeta_\pair}=\varphi_x\one_{\zeta_\pair}
\text{ for all }\al_i\in\smallPi_{\beta(x)}.$
Thus
we have an $\H$-homomorphism
\begin{equation}
  \label{eq;mapx}
  \bar\Y_{\beta(x)}(x\overline{(\zeta_\pair}))
\to \bar \H \varphi_x\one_{\zeta_\pair}
\end{equation}
defined via
$\one_{x\overline{(\zeta_\pair})}
\mapsto \varphi_x\cdot\one_{\zeta_\pair}$.

It follows from 
 (\ref{eq;con0}) that
 $ \bar\Y_\beta(x)(x\overline{(\zeta_\pair}))
\text{ is irreducible},$
and thus
 $ \bar\Y_\beta(x)(x\overline{(\zeta_\pair}))
\cong \bar \H \varphi_x\one_{\zeta_\pair}.$
Therefore $Q\varphi_x\one_{\zeta_{\pair}}\neq 0$.

Next let us prove that
$\{ Q\varphi_{x} \one_{\zeta_{{\pair}}}\}_{
 x\in X^{\beta_\pair}_\l}$ is
linearly independent.

Note that $\mu\in P_\g^+(\l)$ implies that
$\beta_\pair$ and $\zeta_\pair$ satisfies the condition
(\ref{cond;dominant}) in Proposition \ref{pr;w(zeta)=zeta}.
Hence we have
\begin{equation}
  \label{eq;con5}
  \Y_{\beta_\pair}(\zeta_{\pair})_{\zeta_\pair}^\gen=\C
\one_{\zeta_\pair}.
\end{equation}
Combining with (\ref{eq;con2}), 
it follows that
$$
\text{the map }  X^{\beta_\pair}_\l\to (\t')_\beta^*
\text{ given by }x\mapsto x(\zeta_\pair)
\text{ is injective.}
$$
For any $x\in X^{\beta_\pair}_\l$,
it can be proved from Lemma \ref{weightcon} and Lemma \ref{lem;eta}
that
$$
(x(\zeta_{\pair}),\al_i)\leq 0 
\text{ for any }i=1,\dots,N-1.
$$
Therefore
$$
  \text{the map }
X^{\beta_\pair}_\l\to (\t')_\beta^*/\bar W
\text{ given by } 
x\mapsto [x(\zeta_\pair)]
\text{ is injective.}
$$
Hence 
$\{ Q\varphi_{x}\cdot \one_{\zeta_{{\pair}}}\}_{
 x\in X^{\beta_\pair}_\l}$ is
linearly independent.
\qed

\medskip
The specialized character formula
$(\ref{eq;character})$ and the above
Theorem $\ref{lower}$ imply the following:
\begin{con}\label{con;symmetricpart}
For $\lm,\mu\in\Pg^+(\l)$, we have
  \begin{equation}
    \V (\mu, \lambda)^{\ol W}=
 \bigoplus_{T\in\T_s^{(\l)}({\pair})}
\bigoplus_{\eta\in\ol P_-(w_T)} \C Q\varphi_{t_{\eta}\cdot w_T}
\cdot \one_{\zeta_{{\pair}}}.
\end{equation}
\end{con}

{(Tomoyuki Arakawa)}\ {\sc Graduate School of Mathematics,
Nagoya University, Japan}

e-mail: tarakawa@math.nagoya-u.ac.jp

\smallskip
{(Takeshi Suzuki)}\ {\sc Research Institute for Mathematical Sciences,
Kyoto University, Japan}

e-mail: takeshi@kurims.kyoto-u.ac.jp

\smallskip
{(Akihiro Tsuchiya)}\ {\sc Graduate School of Mathematics,
Nagoya University, Japan}

e-mail: tsuchiya@math.nagoya-u.ac.jp
\end{document}